\documentclass[11pt]{article}
\usepackage{amsfonts,amssymb,amsopn,amsmath,mathrsfs,theorem,accents}
\usepackage{hyperref}
\usepackage{graphicx}
\usepackage{verbatim}
\usepackage{authblk}
\usepackage{color}
\usepackage{enumitem}
\usepackage[normalem]{ulem}
\usepackage{bm}
\usepackage[font=small]{caption}

\usepackage{tikz} 
\usetikzlibrary{math,calc,arrows,positioning,fit,petri,external}

\newcounter{thm_counter}
\numberwithin{thm_counter}{section}

\newtheorem{lemma}[thm_counter]{Lemma}
\newtheorem{prop}[thm_counter]{Proposition}
\newtheorem{theorem}[thm_counter]{Theorem}

\newtheorem{cor}[thm_counter]{Corollary}

\theorembodyfont{\upshape}
\newtheorem{remark}[thm_counter]{Remark}

\numberwithin{equation}{section} 
\numberwithin{figure}{section} 

\newcommand{\dis}{\displaystyle}

\newcommand{\noi}{\noindent}
\newcommand{\halmos}{\rule{1ex}{1.4ex}}
\newcommand{\QED}{\nopagebreak{\hspace*{\fill}$\halmos$\medskip}}

\newcommand{\quand}{\quad\mbox{and}\quad}

\newcommand{\bt}{\begin{theorem}}
\newcommand{\et}{\end{theorem}}
\newcommand{\bl}{\begin{lemma}}
\newcommand{\el}{\end{lemma}}
\newcommand{\bp}{\begin{prop}}
\newcommand{\ep}{\end{prop}}
\newcommand{\bcor}{\begin{cor}}
\newcommand{\ecor}{\end{cor}}
\newcommand{\br}{\begin{remark}}
\newcommand{\er}{\end{remark}}
\newcommand{\bcon}{\begin{conjecture}}
\newcommand{\econ}{\end{conjecture}}

\newenvironment{Proof}[1][]{\noi\textbf{Proof #1}}{\QED}
\newcommand{\bpro}{\begin{Proof}}
\newcommand{\epro}{\end{Proof}}

\newcommand{\be}{\begin{equation}}
\newcommand{\ee}{\end{equation}}
\newcommand{\ba}{\begin{array}}
\newcommand{\ea}{\end{array}}
\newcommand{\bc}{\be\begin{array}{r@{\,}c@{\,}l}}
\newcommand{\bac}{\begin{array}{r@{\,}c@{\,}l}}
\newcommand{\ec}{\end{array}\ee}

\newcommand{\ga}{\gamma}
\newcommand{\Ga}{\Gamma}
\newcommand{\de}{\delta}
\newcommand{\De}{\Delta}
\newcommand{\eps}{\varepsilon}

\newcommand{\sig}{\sigma}

\newcommand{\Ai}{{\cal A}}

\newcommand{\Ci}{{\cal C}}

\newcommand{\Ki}{{\cal K}}

\newcommand{\Mi}{{\cal M}}

\newcommand{\Si}{{\cal S}}
\newcommand{\Ti}{{\cal T}}

\newcommand{\Xc}{{\cal X}}
\newcommand{\Yi}{{\cal Y}}

\newcommand{\li}{\langle}
\newcommand{\re}{\rangle}

\newcommand{\up}{\uparrow}
\newcommand{\down}{\downarrow}
\newcommand{\updo}{\updownarrow}
\newcommand{\sub}{\subset}
\newcommand{\beh}{\backslash}

\newcommand{\ov}{\overline}

\newcommand{\pa}{\partial}

\newcommand{\cn}{\colon}
\newcommand{\di}{\mathrm{d}}
\newcommand{\half}{{[0,\infty)}}

\def\to{\rightarrow} 

\def\sw{\subseteq} 
\def\mc{\mathcal} 
\def\mb{\mathbb} 
\def\sc{\setminus} 
\def\P{\mb{P}}
\def\Sp{\mb{S}}
\def\R{\mb{R}} 

\def\N{\mb{N}}

\def\Z{\mb{Z}}

\def\~{\sim}
\def\-{\,;\,} 

\def\li{\langle}

\def\qed{$\blacksquare$}
\def\1{1}
\def\cadlag{cadlag}

\def\l{\left}
\def\r{\right}

\fontsize{11pt}{16.0pt}
\selectfont

\def\epsilon{\varepsilon} 

\newcommand{\pre}{\preceq}

\newcommand{\dr}{d_{\ov\R}}
\newcommand{\Cb}{{\mathbf C}}
\newcommand{\lef}{\vartriangleleft}
\newcommand{\mfs}{\mathfrak{s}}
\newcommand{\Cl}{{\rm Cl}}

\begin{document}

\hyphenation{equip-ped}
\hyphenation{re-de-rive}

\makeatletter\@addtoreset{equation}{section}
\makeatother\def\theequation{\thesection.\arabic{equation}} 

\renewcommand{\labelenumi}{{\rm (\roman{enumi})}}
\renewcommand{\theenumi}{\roman{enumi}}

\allowdisplaybreaks

\author[1]{Nic Freeman\thanks{University of Sheffield, n.p.freeman@sheffield.ac.uk, https://nicfreeman1209.github.io/Website/}}
\author[2]{Jan M.\ Swart\thanks{Institute of Information Theory and Automation, Pod vod\'arenskou v\v{e}\v{z}\'i 4, 18200 Praha 8, Czech Republic, swart@utia.cas.cz}}
\affil[1]{University of Sheffield}
\affil[2]{The Czech Academy of Sciences, Institute of Information Theory and Automation}

\title{Tightness criteria for random compact sets of {\cadlag} paths}
\date{\today}
\maketitle

\begin{abstract}
We give tightness criteria for random variables taking values in the space of all compact sets of {\cadlag} real-valued paths, in terms of both the Skorokhod J1 and M1 topologies. This extends earlier work motivated by the study of the Brownian web that was concerned only with continuous paths. In the M1 case, we give a natural extension of our tightness criteria  which ensures that non-crossing systems of paths have weak limit points that are also non-crossing. This last result is exemplified through a rescaling of heavy-tailed Poisson trees
and a more general application to weaves.
\end{abstract}

\vspace{.5cm}

\noi
{\it MSC 2020.} Primary: 60B10; Secondary: 60G07.\\
%
{\it Keywords.} Tightness, Skorokhod topology, J1 topology, M1 topology, path space, non-crossing paths.\\
{\it Acknowledgments.} Work sponsored by GA\v{C}R grant 22-12790S.

\tableofcontents


\section{Main results}

\subsection{Introduction}\label{S:intro}

A central theme in probability theory is the weak convergence in law of stochastic processes, the canonical example of which is the rescaling of random walks to Brownian motion and, more generally, to $\alpha$-stable processes. The modern perspective is to treat a stochastic process as a single random variable, a random path, and to consider convergence in terms of probability laws on the space of all possible paths of the process. Paths are usually assumed to be \emph{{\cadlag}}, that is, right-continuous with left limits, with time domain $[0,\infty)$ and values in some Polish space $M$. The space $D_{[0,\infty)}(M)$ of such paths is often equipped with Skorokhod's J1 topology, which is commonly known as `the' Skorokhod topology. For real-valued processes, the slightly weaker Skorokhod M1 topology is sometimes required instead; see, e.g., \cite{Whi02,DGM24}. Proofs of weak convergence typically involve establishing that (1) a sequence of probability laws on the space of all paths is tight and (2) has a unique weak cluster point; the classical works \cite{EK86} and \cite{Bil99} are standard references for this approach. For this reason, tightness criteria are of central importance.

The introduction of the Brownian web \cite{FINR04} started interest in random variables that are not a single random path, but rather a random compact \emph{set of} paths. The framework for doing so established by \cite{FINR04} for the Brownian web (often known as `the Brownian web topology') treats continuous paths. More recently random compact sets whose elements are paths with jumps have been considered \cite{MRV19,FS24}, along with surprising tightness failures for continuous paths that are related to the appearance of a small fraction of jumps in the limit \cite{BMRV06, SSY21}. This naturally asks for tightness criteria for random sets of {\cadlag} paths, which is the focus of the present article.
 

We consider systems in which different paths are typically defined on different time domains. The Brownian web, for example, covers space-time $\R\times\R$ with continuous paths $\pi$ for which the time domains are of the form $[s,\infty)$, and $s\in\R$ is called the starting time of $\pi$. This necessitates an underlying topology (on individual paths) that allows for the convergence of a sequence of paths that may all have different time domains. In the work of \cite{FINR04}, for continuous paths, such convergence can roughly be described as locally uniform convergence of functions plus convergence of the starting times. In \cite{FS25} we generalised this principle by introducing a space of {\cadlag} paths for which the time domains are \emph{arbitrary} closed subsets of the real line. Although our focus is on random sets of paths, we note that this generalization is natural in even the most basic examples, for example, we may then express the rescaling of discrete time random walks to continuous time processes without involving (linear or otherwise) interpolation in between discrete times of the walk.

In \cite{FS25} we showed that the space introduced therein of {\cadlag} paths (with arbitrary closed subsets as time domains) was Polish, for both the J1 and M1 cases.
Our present paper builds on this work by deriving tightness criteria for random compact subsets of this space, in the case where paths are real-valued. The main result is Theorem~\ref{T:tight} below. It provides a robust criterion that involves checking a single property, based on how frequently paths cross small intervals of $\R$.

When working with {\cadlag} paths on varying time domains, one complication is that the property of two paths not crossing each other is not automatically preserved by taking limits. This complication arises from jumps that may form at the starting times of limiting paths; we remark that our state space permits this possibility, because it is necessarily a feature of treating random sets of {\cadlag} paths as compact sets (see \cite[Figure 1.1.1]{FS24} or the discussion in \cite[Appendix A.5]{FS24} for details on this point). The theory that has developed around the Brownian web has shown that systems of non-crossing paths are substantially easier to work with than systems with crossing paths. With this in mind, in Theorem~\ref{T:noncros} we note a natural analogue of our main result which implies not only tightness, but also that all weak limit points are non-crossing.

In Section~\ref{s:application} we give two applications of our results. We start by giving tightness criteria for weaves -- a \emph{weave}, introduced in \cite{FS24}, is a random compact subset $\mc{A}$ consisting of paths $\pi\in\mc{A}$ that do not cross one another, and for which, with probability one, for all $z\in\R^2$ there exists $\pi\in\mc{A}$ such that $z$ is contained within the graph of $\pi$; that is, the paths cover space-time. We show that, for weaves, tightness essentially comes down to tightness of the motion of a single particle. In Section \ref{s:application_poisson_trees} we exemplify this result on a sequence of rescaled Poisson trees in which the limiting motion of a single path is an $\alpha$-stable process.

Our main results are proved in Sections \ref{S:tightproof} and \ref{S:crosproof}. Section~\ref{S:tightproof} contains the proof of Theorems \ref{T:ctight} and \ref{T:tight}, and Section~\ref{S:crosproof} contains the proof of Theorem~\ref{T:noncros}. The results of Section~\ref{s:application} are proved on the spot, to demonstrate the application of our main results.


\subsection{Path space}\label{S:topol}

We begin by recalling the necessary definitions and related background from \cite{FS25}. Let $\ov\R:=[-\infty,\infty]$. By \cite[Lemma~2.6]{FS25}, there exists a unique metrisable topology on
\be
\R^2_{\rm c}:=(\ov\R\times\R)\cup\big\{(\ast,-\infty),(\ast,+\infty)\big\}
\ee
such that $\R^2_{\rm c}$ is compact and a sequence $(x_n,t_n)\in\R^2_{\rm c}$ converges to a limit $(x,t)$ if and only if
\begin{enumerate}
\item $t_n\to t$ in the topology on $\ov\R$,
\item if $t\in\R$, then $x_n\to x$ in the topology on $\ov\R$.
\end{enumerate}
The space $\R^2_{\rm c}$ has earlier been introduced in \cite{FINR04}. Keeping with tradition, we think of the first coordinate as space and the second as time, so in pictures time runs upwards. We can think of $\R^2_{\rm c}$ as being obtained from the space $\ov\R^2$ after squeezing the sets $\ov\R\times\{\pm\infty\}$ into the single points $(\ast,\pm\infty)$. For this reason, we call $\R^2_{\rm c}$ the \emph{squeezed space}. We let $d_{\rm sqz}$ denote any metric generating the topology on~$\R^2_{\rm c}$. The metric $d_{\rm sqz}$ cares less about the spatial distance between points if their time coordinates are large, which in what follows will mean that two paths are close if their graphs are locally close.

Let $\hat I$ be a nonempty closed subset of $\ov\R$, let $I:=\hat I\cap\R$, and let
\be
I^-:=\big\{t\in I:(t-\eps,t)\cap I\neq\emptyset\ \forall\eps>0\big\}
\quand
I^+:=\big\{t\in I:(t,t+\eps)\cap I\neq\emptyset\ \forall\eps>0\big\}
\ee
denote the sets of points in $I$ that can be approximated from the right or left, respectively. Let $\pi^\pm\cn I\to\ov\R$ be functions such that $\pi^-$ is left-continuous, $\pi^+$ is right-continuous, and
\be\label{pipm}
\pi^-(t)=\lim_{s\up t}\pi^+(s)\quad(t\in I^-)\quand\pi^+(t)=\lim_{s\down t}\pi^-(s)\quad(t\in I^+).
\ee
We call the triple $(\hat I,\pi^-,\pi^+)$ a \emph{path} and we say that $I$ and $\hat I$ are its \emph{domain} and \emph{extended domain}. Roughly, we can think of $\pi^+$ as a {\cadlag} function and of $\pi^-$ as its left-continuous modification, but contrary to the usual conventions we allow $\pi^-$ to differ from $\pi^+$ in points that cannot be approximated from the left.

For $x,z\in\ov\R$ we write $[x,z]:=\{y\in\ov\R:x\wedge z\leq y\leq x\vee z\}$. For a given path $(\hat I,\pi^-,\pi^+)$, we set
\bc
\label{piset}
\dis\pi&:=&\dis\big\{(x,t):t\in I,\ x\in\{\pi(t-),\pi(t+)\}\big\}\cup\big\{(\ast,t):t\in\hat I\beh I\big\},\\[5pt]
\dis\ov\pi&:=&\dis\big\{(x,t):t\in I,\ x\in[\pi(t-),\pi(t+)]\big\}\cup\big\{(\ast,t):t\in\hat I\beh I\big\}.
\ec
We call $\pi$ the \emph{closed graph} and $\ov\pi$ the \emph{filled graph} of $(\hat I,\pi^-,\pi^+)$. We equip these sets with a total order $\pre$ such that two elements $(x,s),(y,t)\in\pi,\ov\pi$ are strictly ordered $(x,s)\prec(y,t)$ if and only if either $s<t$, or $s=t\in\R$ and $x$ lies closer to $\pi(t-)$ than $y$.

It is possible to give an alternative definition of paths directly in terms of their closed or filled graphs. Let $\Pi$ denote the space of all nonempty compact subsets $\pi$ of $\R^2_{\rm c}$ that are equipped with a total order $\pre$ such that
\begin{enumerate}
\item $\pi^{\li 2\re}:=\big\{(z,z')\in\pi^2:z\pre z'\big\}$ is a closed subset of $\pi^2:=\pi\times\pi$, where we equip $\pi$ with the induced topology from$\R^2_{\rm c}$ and $\pi^2$ with the product topology.
\item If $(x,s),(x',s')\in\pi$ satisfy $s<s'$, then $(x,s)\pre(x',s')$.
\item For each $t\in\R$, the set $\{x\in\Xc:(x,t)\in\pi\}$ has at most two elements.
\end{enumerate}
We define $\ov\Pi$ similarly, but with (iii) replaced by
\begin{itemize}
\item[(iii)'] For each $t\in\R$, the set $\{x\in\Xc:(x,t)\in\pi\}$ is an interval.
\end{itemize}
The following lemma has been proved in \cite[Lemma~1.2]{FS25}.

\bl[Closed and filled graphs]
The\label{L:graph} closed graph $\pi$ and filled graph $\ov\pi$ of a path $(\hat I,\pi^-,\pi^+)$ satisfy $\pi\in\Pi$ and $\ov\pi\in\ov\Pi$. Conversely, each element of $\Pi$ is the closed graph of a unique path $(\hat I,\pi^-,\pi^+)$, and each element of $\ov\Pi$ is the filled graph of a unique path $(\hat I,\pi^-,\pi^+)$.
\el

In view of Lemma~\ref{L:graph}, we usually identify a path $(\hat I,\pi^-,\pi^+)$ with its closed graph $\pi$ and let $\ov\pi$ denote its associated filled graph. We let $I_\pi$ and $\hat I_\pi$ denote the domain and extended domain of $\pi$. Instead of $\pi^\pm(t)$ we write $\pi(t\pm)$ which is suggestive of left- and right limits, as is justified by (\ref{pipm}), at least in points that can be approximated from the left or right. We also set $\pi(\pm\infty):=\ast$ if $\pm\infty\in\hat I$.

Following \cite{FS25}, we will equip $\Pi$ and $\ov\Pi$ with topologies that are generalisations of Skorokhod's J1 and M1 topologies. A \emph{correspondence} between two sets $A,B$ is a set $C\sub A\times B$ such that
\be\label{corresp}
\forall a\in A\ \exists b\in B\mbox{ s.t.\ }(a,b)\in C
\quand
\forall b\in B\ \exists a\in A\mbox{ s.t.\ }(a,b)\in C.
\ee
We let ${\rm Corr}(A,B)$ denote the set of all correspondences between $A$ and $B$. If $\Xc$ is a metrisable topological space, then we let $\Ki_+(\Xc)$ denote the space of nonempty compact subsets of $\Xc$. If $d$ is a metric generating the topology on $\Xc$, then the corresponding \emph{Hausdorff metric} on $\Ki_+(\Xc)$ is defined as
\be
d_{\rm H}(K_1,K_2):=\inf_{z_1\in K_1}d(z_1,K_2)\vee\inf_{z_2\in K_2}d(z_2,K_1)
=\inf_{C\in{\rm Corr}(K_1,K_2)}\sup_{(z_1,z_2)\in C}d(z_1,z_2),
\ee
where as usual $d(z,K):=\inf_{z'\in K}d(z,z')$ denotes the distance of a point to a set. It follows from \cite[Lemma~B.1]{SSS14} that the topology on $\Ki_+(\Xc)$ only depends on the topology on $\Xc$ and not on the choice of the metric. We call this the \emph{Hausdorff topology} on $\Ki_+(\Xc)$. We note that \cite[Lemmas B.2 and B3]{SSS14} show that $\Ki_+(\Xc)$ is Polish if $\Xc$ is Polish, and compact if $\Xc$ is compact.

For paths $\pi_1,\pi_2\in\Pi$, we let ${\rm Corr}_+(\pi_1,\pi_2)$ denote the set of correspondences $C$ between $\pi_1$ and $\pi_2$ that are \emph{monotone} in the sense that
\be
\mbox{there are no $(z_1,z_2),(z'_1,z'_2)\in C$ such that $z_1\prec_1 z'_1$ and $z'_2\prec_2 z_2$}.
\ee
If $d_{\rm sqz}$ is a metric generating the topology on the squeezed space $\R^2_{\rm c}$, then we define a corresponding metric $d_{\rm J1}$ on the path space $\Pi$ by
\be\label{J1def}
d_{\rm J1}(\pi_1,\pi_2):=\inf_{C\in{\rm Corr}_+(\pi_1,\pi_2)}\sup_{(z_1,z_2)\in C}d_{\rm sqz}(z_1,z_2)\qquad\big(\pi_1,\pi_2\in\Pi\big).
\ee
We define $d_{\rm M1}(\pi_1,\pi_2)$ in the same way, but with ${\rm Corr}_+(\pi_1,\pi_2)$ replaced by ${\rm Corr}_+(\ov\pi_1,\ov\pi_2)$, the set of monotone correspondences between the filled graphs $\ov\pi_1$ and $\ov\pi_2$. It follows from \cite[Thm~1.9]{FS25} that the topologies generated by these metrics do not depend on the choice of the metric $d_{\rm sqz}$ on $\R^2_{\rm c}$, and the space $\Pi$, equipped with either topology, is Polish.\footnote{We adhere to the nowadays fairly standard convention that a topological space is Polish if it is separable and there exists a complete metric generating the topology. For the space $\Pi$, the fact that $d_{\rm J1}$ and $d_{\rm M1}$ generate Polish topologies does not imply that these metrics themselves are complete, which, indeed, they are not.} By \cite[Thms 1.4 and 1.9]{FS25}, when restricted to the classical spaces of {\cadlag} paths that are defined on a fixed domain, the metric $d_{\rm J1}$ generates Skorokhod's J1-topology while $d_{\rm M1}$ generates Skorokhod's M1-topology.

We let $\Pi_{\rm c}$ and $\Pi^|$ denote the subspaces of $\Pi$ consisting of \emph{continuous} and \emph{connected} paths, respectively, defined as
\bc\label{Pic}
\dis\Pi_{\rm c}&:=&\dis\big\{\pi\in\Pi:\pi(t-)=\pi(t+)\ \forall t\in I_\pi\big\},\\[5pt]
\dis\Pi^|&:=&\dis\big\{\pi\in\Pi:\hat I_\pi\mbox{ is an interval}\big\}.
\ec
We also write $\Pi^\up$ and $\Pi^\down$ for the subsets of $\Pi^|$ consisting of paths whose domain is an interval that is unbounded from above and below, respectively, and we let $\Pi^\updo:=\Pi^\up\cap\Pi^\down$ denote the set of bi-infinite paths. We naturally write $\Pi^|_{\rm c}:=\Pi_{\rm c}\cap\Pi^|$ and so on. Then $\Pi^\up_{\rm c}$ is the classical path space of the Brownian web \cite{FINR04}. For paths $\pi\in\Pi_{\rm c}$ we write $\pi(t):=\pi(t-)=\pi(t+)$ $(t\in I_\pi)$. Since for continuous paths, the closed and filled graphs coincide, the metrics $d_{\rm J1}$ and $d_{\rm M1}$ generate the same topology on $\Pi_{\rm c}$. It follows from \cite[Thm~1.4]{FS25} that $\Pi^|_{\rm c}$ is closed as a subset of $\Pi$ in the J1-topology. We note that paths in $\Pi^\up$ can make a jump at their starting time. This is different from the usual conventions for the space of {\cadlag} functions defined on $\half$ but will be crucial for Theorem~\ref{T:tight} below as well as in future applications we have in mind.

\subsection{Tightness for random compact sets of paths}\label{S:tight}

Let $\Pi$ be the path space introduced in the previous subsection, equipped with either the J1 or M1 topology, and let $\Ki_+(\Pi)$ be the space of nonempty compact subsets of $\Pi$, equipped with the Hausdorff topology. We note that since $\Pi$ is Polish under both the J1 and M1 topologies, by \cite[Lemma~2.7]{FS25} the resulting topology on $\Ki_+(\Pi)$ is also Polish in either case. We will be interested in weak convergence of probability measures on $\Ki_+(\Pi)$ with respect to both the J1 and M1 topologies on $\Pi$. General topology tells us\footnote{In any metrisable space, a sequence $x_n$ converges to a limit $x$ if and only if the set $\{x_n:n\in\N\}$ is precompact (i.e., its closure is compact) and $x$ is its only cluster point (i.e., subsequential limit point).} that a sequence of probability measures $\mu_n$ on $\Ki_+(\Pi)$ converges weakly to a limit if and only if the set $\{\mu_n:n\in\N\}$ is precompact and $\mu$ is its only cluster point. Since $\Ki_+(\Pi)$ is Polish, Prokhorov's theorem tells us that a family $\{\mu_\ga:\ga\in\Ga\}$ of probability measures on $\Ki_+(\Pi)$ is precompact if and only if it is \emph{tight}, that is, for all $\eps>0$, there exists a compact $\Ci\sub\Ki_+(\Pi)$ such that
\be\label{tightdef}
\sup_{\ga\in\Ga}\mu_\ga\big(\Ki_+(\Pi)\beh\Ci\big)\leq\eps.
\ee
We are therefore naturally interested in tightness criteria for families of probability measures on $\Ki_+(\Pi)$, with respect to the J1 and M1 topologies on $\Pi$. As a warm-up, we discuss tightness on $\Ki_+(\Pi_{\rm c})$, where $\Pi_{\rm c}$, defined in (\ref{Pic}), is the space of continuous paths. As already mentioned, the metrics $d_{\rm J1}$ and $d_{\rm M1}$ generate the same topology on $\Pi_{\rm c}$.

For each $\pi\in\Pi$ and real $T,\de>0$, we set
\bc\label{De2}
\dis\De^2_{T,\de}(\pi)&:=&\dis\big\{(x,y)\in\ov\R^2:\exists-T\leq s\leq t\leq T\mbox{ s.t.\ }t-s\leq\de\\
&&\dis\phantom{\big\{(x,y):}\mbox{and }(x,s),(y,t)\in\pi,\ (x,s)\pre(y,t)\big\}.\ec
For each $T,\de,\eps>0$ and $r\in\R$, we define sets of paths by
\bc
\dis\Si^+_{T,\de,\eps,r}&:=&\big\{\pi\in\Pi:\exists(x,y)\in\De^2_{T,\de}(\pi)\mbox{ s.t.\ }x\leq r,\ r+\eps\leq y\big\},\\[5pt]
\dis\Si^-_{T,\de,\eps,r}&:=&\big\{\pi\in\Pi:\exists(x,y)\in\De^2_{T,\de}(\pi)\mbox{ s.t.\ }y\leq r,\ r+\eps\leq x\big\},
\ec
and we set $\Si^2_{T,\de,\eps,r}:=\Si^+_{T,\de,\eps,r}\cup\Si^-_{T,\de,\eps,r}$. 
The set $\dis\Si^+_{T,\de,\eps,r}$ comprises paths that cross $[r,r+\eps]$
from left to right within time $[-T,T]$,
with the condition that once the crossing has begun it takes time less than $\delta$ to complete.
The set $\dis\Si^-_{T,\de,\eps,r}$ corresponds to crossings from right to left.
Our first result is an extension of earlier work \cite[Prop~B1]{FINR04} which therein is restricted to the space $\Pi^\up_{\rm c}$.

\bt[Tightness criterion for sets of continuous paths]
Let\label{T:ctight} $(\Ai_\ga)_{\ga\in\Ga}$ be a family of random variables with values in $\Ki_+(\Pi_{\rm c})$. Then the laws $(\mu_\ga)_{\ga\in\Ga}$ with $\mu_\ga:=\P[\Ai_\ga\in\,\cdot\,]$ are tight with respect to the topology on $\Pi_{\rm c}$ if and only if
\be\label{ctight}
\lim_{\de\to 0}\sup_{\ga\in\Ga}\P\big[\Si^2_{T,\de,\eps,r}\cap\Ai_\ga\neq\emptyset\big]=0\quad\forall T,\eps>0,\ r\in\R.
\ee
\et

We now turn our attention to paths with jumps. For each $\pi\in\Pi$ and real $T,\de>0$, we set
\bc\label{De3}
\dis\De^3_{T,\de}(\pi)&:=&\dis\big\{(x,y,z)\in\ov\R^3:\exists-T\leq s\leq t\leq u\leq T\mbox{ s.t.\ }u-s\leq\de\\
&&\dis\phantom{\big\{(x,y,z):}\mbox{and }(x,s),(y,t),(z,u)\in\pi,\ (x,s)\pre(y,t)\pre(z,u)\big\}.
\ec
For each $T,\de,\eps>0$ and $r\in\R$, we define sets of paths by
\bc\label{Sipm}
\dis\Si^{+-}_{T,\de,\eps,r}&:=&\big\{\pi\in\Pi:\exists(x,y,z)\in\De^3_{T,\de}(\pi)\mbox{ s.t.\ }x,z\leq r,\ r+\eps\leq y\big\},\\[5pt]
\dis\Si^{-+}_{T,\de,\eps,r}&:=&\big\{\pi\in\Pi:\exists(x,y,z)\in\De^3_{T,\de}(\pi)\mbox{ s.t.\ }y\leq r,\ r+\eps\leq x,z\big\},\\[5pt]
\dis\Si^{++}_{T,\de,\eps,r}&:=&\big\{\pi\in\Pi:\exists(x,y,z)\in\De^3_{T,\de}(\pi)\mbox{ s.t.\ }x\leq r,\ r+\eps\leq y\leq r+2\eps,\ r+3\eps\leq z\big\},\\[5pt]
\dis\Si^{--}_{T,\de,\eps,r}&:=&\big\{\pi\in\Pi:\exists(x,y,z)\in\De^3_{T,\de}(\pi)\mbox{ s.t.\ }z\leq r,\ r+\eps\leq y\leq r+2\eps,\ r+3\eps\leq x\big\},
\ec
and we set
\be
\Si^{\rm J}_{T,\de,\eps,r}:=\Si^{++}_{T,\de,\eps,r}\cup\Si^{--}_{T,\de,\eps,r}
\quand
\Si^{\rm M}_{T,\de,\eps,r}:=\Si^{+-}_{T,\de,\eps,r}\cup\Si^{-+}_{T,\de,\eps,r}.
\ee
The set $\dis\Si^{+-}_{T,\de,\eps,r}$ comprises paths that cross $[r,r+\eps]$
from left to right, and then back again from right to left, during time $[-T,T]$
with the whole exercise taking time less that $\de$;
whilst $\dis\Si^{-+}_{T,\de,\eps,r}$ corresponds to similar movement in opposite directions.
The set $\dis\Si^{++}_{T,\de,\eps,r}$
comprises paths that cross $[r,r+\eps]$ from left to right,
then take a value in $[r+\eps,r+2\eps]$,
then cross $[r+2\eps,r+3\eps]$ from left to right, during time $[-T,T]$
with the whole exercise taking time less that $\de$;
whilst $\dis\Si^{--}_{T,\de,\eps,r}$ corresponds to similar movement in opposite directions.

Our main result is the following theorem. Roughly, this says that the laws of $(\Ai_\ga)_{\ga\in\Ga}$ are tight with respect to the J1 topology on $\Pi$ if the sets $\Ai_\ga$ do not contain paths that make two jumps in an arbitrarily brief time after each other. For the M1 topology, it suffices to look only at two jumps in opposite directions.

\bt[Tightness criteria for sets of {\cadlag} paths]\hspace{5pt}
Let\label{T:tight} $(\Ai_\ga)_{\ga\in\Ga}$ be a family of random variables with values in $\Ki_+(\Pi)$. Then the laws $(\mu_\ga)_{\ga\in\Ga}$ with $\mu_\ga:=\P[\Ai_\ga\in\,\cdot\,]$ are tight with respect to the J1 topology on $\Pi$ if and only if
\be\ba{r@{\ }l}\label{tight}
{\rm(i)}&\dis\lim_{\de\to 0}\sup_{\ga\in\Ga}\P\big[\Si^{\rm M}_{T,\de,\eps,r}\cap\Ai_\ga\neq\emptyset\big]=0\quad\forall T,\eps>0,\ r\in\R,\\[5pt]
{\rm(ii)}&\dis\lim_{\de\to 0}\sup_{\ga\in\Ga}\P\big[\Si^{\rm J}_{T,\de,\eps,r}\cap\Ai_\ga\neq\emptyset\big]=0\quad\forall T,\eps>0,\ r\in\R.
\ec
For tightness with respect to the M1 topology on $\Pi$ it is necessary and sufficient that only (\ref{tight}) (i) holds.
\et

Note that if the laws of the $\Ai_\ga$'s are invariant under translations, then it suffices to check (\ref{tight}) (i) and (ii) for $r=0$. If the laws are invariant under reflection then we can replace $\Si^{\rm M}_{T,\de,\eps,r}$ by $\Si^{+-}_{T,\de,\eps,r}$ and $\Si^{\rm J}_{T,\de,\eps,r}$ by $\Si^{++}_{T,\de,\eps,r}$.

\subsection{Non-crossing sets of paths}\label{S:noncros}

Recall from (\ref{Pic}) that $\Pi^|$ is the space of `connected' paths, whose domain is an interval, and $\Pi^\updo$ is the space of bi-infinite paths. Note that $\pi\in\Pi^|$ if and only if its filled graph $\ov\pi$ as defined in (\ref{piset}) is connected in the topological sense. We say that a path $\pi'$ \emph{extends} a path $\pi$ if $\ov\pi\sub\ov\pi'$. Following \cite{FS24}, for $\pi_1,\pi_2\in\Pi^|$, we write $\pi_1\lef\pi_2$ if $\pi_1$ and $\pi_2$ can be extended to bi-infinite paths $\pi'_1,\pi'_2\in\Pi^\updo$ such that $\pi'_1(t\pm)\leq\pi'_2(t\pm)$ for all $t\in\R$. We note that in spite of the suggestive notation, this relation is not transitive, that is, $\pi_1\lef\pi_2\lef\pi_3$ does not imply $\pi_1\lef\pi_3$.\footnote{For example, $\pi_2$ may be the trivial path that contains only the point $(\ast,-\infty)$, in which case $\pi_1\lef\pi_2\lef\pi_3$ holds trivially for arbitrary $\pi_1$ and $\pi_3$.} We say that a set $\Ai\sub\Pi^|$ is \emph{non-crossing} if each $\pi_1,\pi_2\in\Ai$ satisfy $\pi_1\lef\pi_2$ or $\pi_2\lef\pi_1$ (or both). Throughout the present subsection, we equip $\Pi^|$ with the M1 topology and we equip
\be
\Ki_{\rm nc}(\Pi^|):=\big\{\Ai\in\Ki_+(\Pi^|):\Ai\mbox{ is non-crossing}\big\}
\ee
with the corresponding Hausdorff topology. We will prove a tightness criterion for random variables with values in $\Ki_{\rm nc}(\Pi^|)$. We note that when we say that a collection of probability measures is \emph{tight on} $\Ki_{\rm nc}(\Pi^|)$, we mean that the compact subsets that occur in the definition of tightness (see (\ref{tightdef})) are subsets of $\Ki_{\rm nc}(\Pi^|)$ (and not of some larger space). In particular, this implies that each weak cluster point is concentrated on $\Ki_{\rm nc}(\Pi^|)$, which would not follow from the tightness criterion with respect to the M1 topology of Theorem~\ref{T:tight}.

Nevertheless, the tightness criterion for non-crossing sets of paths turns out to be very similar to condition (\ref{tight})~(i), though a bit stronger. To formulate it, for each $\pi_1,\pi_2\in\Pi$ and real $T,\de>0$, we set
\bc\label{De22}
\dis\De^2_{T,\de}(\pi_1,\pi_2)&:=&\dis\big\{(x_1,y_1,x_2,y_2)\in\ov\R^4:\exists-T\leq s_i\leq t_i\leq T\mbox{ s.t.\ }(x_i,s_i),(y_i,t_i)\in\pi_i,\\
&&\dis\hspace{2cm}\ (x_i,s_i)\pre(y_i,t_i)\ (i=1,2)\mbox{ and }(t_1\vee t_2)-(s_1\wedge s_2)\leq\de\big\},
\ec
and we define
\be
\dis\Ci^{\rm M}_{T,\de,\eps,r}:=\big\{(\pi_1,\pi_2)\in\Pi^2:\exists(x_1,y_1,x_2,y_2)\in\De^2_{T,\de}(\pi_1,\pi_2)\mbox{ s.t.\ }x_1,y_2\leq r,\ r+\eps\leq y_1,x_2\big\}.
\ee

\bt[Tightness criterion for non-crossing sets of paths]
Let\label{T:noncros} the space $\Pi^|$ be equipped with the M1 topology and $\Ki_{\rm nc}(\Pi^|)$ with the corresponding Hausdorff topology. Let $(\Ai_\ga)_{\ga\in\Ga}$ be a family of random variables with values in $\Ki_{\rm nc}(\Pi^|)$. Then the laws $(\mu_\ga)_{\ga\in\Ga}$ with $\mu_\ga:=\P[\Ai_\ga\in\,\cdot\,]$ are tight on $\Ki_{\rm nc}(\Pi^|)$ if and only if
\be\label{crotig}
\lim_{\de\to 0}\sup_{\ga\in\Ga}\P\big[\Ci^{\rm M}_{T,\de,\eps,r}\cap(\Ai_\ga\times\Ai_\ga)\neq\emptyset\big]=0\quad\forall T,\eps>0,\ r\in\R.
\ee
\et
In words, the event $\Ci^{\rm M}_{T,\de,\eps,r}\cap(\Ai_\ga\times\Ai_\ga)\neq\emptyset$ occurs when $\mc{A}_\gamma$ contains a path that crosses $[r,r+\eps]$ in one direction, and another (perhaps different) path that crosses $[r,r+\eps]$ in the opposite direction, with both events taking place within a single time window of length $\de$ between times $-T$ and $T$. As already discussed at the end of Section \ref{S:intro}, in Section \ref{s:application} we give some applications of Theorem \ref{T:noncros}.

It is easy to see that (\ref{crotig}) implies (\ref{tight})~(i) (see Lemma~\ref{L:CS} below), but the converse implication does not hold in general. For sets of bi-infinite paths, the two conditions are equivalent. This is a consequence of Theorems \ref{T:tight} and \ref{T:noncros} and the following lemma.

\bl[Sets of bi-infinite paths]
Let\label{L:biclos} $\Pi$ be equipped with the M1 topology and $\Ki_+(\Pi)$ with the corresponding Hausdorff topology. Then $\Ki_{\rm nc}(\Pi^\updo)$ is a closed subset of $\Ki_+(\Pi)$.
\el

\section{Applications}
\label{s:application}

Throughout Section \ref{s:application} we will use a piece of notation that is in common usage for the Brownian web, and related objects,
but that we have not so far introduced.
For $\pi\in\Pi^\uparrow$ we write $\sigma_\pi$ for the initial time of $\pi$, that is, if $\pi$ has time domain $[t,\infty)$ then $\sigma_\pi=t$.

\subsection{Weaves}
\label{s:application_noncr}

In this section we prove a result that specialises Theorem \ref{T:noncros}.
A \emph{weave}, introduced and studied in \cite{FS24},
is a random compact subset $\mc{A}\sw\Pi^\uparrow$
that is non-crossing and for which,
with probability one, for all $z\in\R^2$ there exists $\pi\in\mc{A}$ such that $z\in\ov{\pi}$, that is, the paths cover space-time.
We consider weaves with a law that is invariant under deterministic translations of space-time.
In this setting, we show that tightness of the motion of a single particle in the sense of 
the \emph{classical} Skorokhod space $D_{[0,\infty)}(\R)$ 
(as defined in e.g.~\cite{EK86, Bil99, Whi02}) 
essentially controls tightness of the whole system in $\mc{K}(\Pi^\uparrow)$. 
We show also that translation invariance can be replaced by a more technical condition that is uniform in space and time.

Before stating the result we must note a subtle technical point.
Let $(X_n)$ be a sequence of $D_{[0,\infty)}(\R)$-valued random variables.
Although $D_{[0,\infty)}(\R)$ is a subspace of $\Pi^\uparrow$,
tightness of $(X_n)$ in $\Pi^\uparrow$ in the sense of Theorem \ref{T:tight}
does not imply tightness of $(X_n)$ in $D_{[0,\infty)}(\R)$ in the classical sense.
The distinction arises because for $\pi\in D_{[0,\infty)}(\R)$ we must have $\pi(0-)=\pi(0+)$
whereas $\Pi$ permits jumps at time $0$,
more precisely
\be\label{Ddef}
D_{[0,\infty)}(\R)=\big\{\pi\in\Pi^\up:\sig_\pi=0,\ \pi(0-)=\pi(0+)\big\}.
\ee
Tightness conditions for $D_{[0,\infty)}(\R)$
(as found in \cite{EK86, Bil99, Whi02} and suchlike)
must therefore enforce that limit points do not jump at time $0$,
whereas tightness conditions (such as Theorems \ref{T:tight} and \ref{T:noncros}) 
for the larger space $\Pi^\uparrow$ need not do so. For single paths the precise connection,
which is a straightforward consequence of \cite[Thm~1.13]{FS25}, is as follows.

\bl[Tightness on $D_{[0,\infty)}(\R)$]
\label{l:D_tightness} 
Let $(X_n)$ be a sequence of random variables with values in $D_{[0,\infty)}(\R)$. Then their laws are tight on $D_{[0,\infty)}(\R)$ under J1 (resp.~M1) if and only if these laws are tight on $\Pi$ under J1 (resp.~M1) and additionally,
for all $\eps>0$ we have
$
\lim_{\de\to 0}\sup_{n\in\N}\P\big[\sup_{t\in[0,\de]}|X_n(t)|\geq\eps\big]=0.
$
\el

In the statement of Lemma \ref{l:D_tightness} we have used the standard notation $X_n(t)$, where $t\in\R$ and $X_n$ is a $D_{[0,\infty)}(\R)$-valued process.
Note in particular that here the time argument is $t$, rather than $t-$ or $t+$.
Hereon, if some $\pi\in\Pi^\up$ starts at time zero and satisfies $\pi(0-)=\pi(0+)$
then we will allow ourselves to view $\pi$ as an element of $D_{[0,\infty)}(\R)$,
writing $\pi(t)=\pi(t+)$ for all $t\geq 0$.
This slight abuse of notation is helpful when applying results concerning the classical Skorohod spaces.

More generally, if $\pi\in\Pi^\up$ is continuous at its initial time $s$, 
then we view $\pi$ as an element of the space $D_{[s,\infty)}(\R)$
of classical {\cadlag} paths with time domain $[s,\infty)$
and write $\pi(t)=\pi(t+)$ for all $t\geq s$.
For a random path $\pi\in D_{[s,\infty)}(\R)$ 
the Markov property with respect to some filtration $(\mc{F}_t)_{t\in\R}$ is defined in the usual way
i.e.~equivalent to the requirement that $t\mapsto\pi(s+t)$,
which is a $D_{[0,\infty)}(\R)$-valued process, 
is a Markov process with respect to $(\mc{F}_{s+t})_{t\geq 0}$.
The strong Markov property is defined similarly.

We now recount some important facts from \cite[Section 2.4]{FS24} concerning weaves.
Let $\mc{A}$ be a weave.
Then, except for $z$ within a (deterministic) Lebesgue null set $R_\mc{A}\sw\R^2$,
for each $z\in\R^2$  with probability one there exists a unique path $\pi_z\in\Pi^\uparrow$ that begins at $z$ and does not cross $\mc{A}$.
We say that $\pi_z$ is the \emph{one-particle motion} of $\mc{A}$ from $z$.

We say that a weave $\mc{A}$ is \emph{homogeneous} if the law of $\mc{A}$ is invariant under deterministic translations of space and time. In this case the set $R_\mc{A}$ is empty:
for each deterministic $z=(x,s)\in\R^2$ there almost surely exists a unique path $\pi_z\in\Pi^\uparrow$ that begins at $z$
and does not cross $\mc{A}$.
The law of $\pi_z$ might include a jump at its starting time, 
but if it does not do so then
\begin{equation}
\label{eq:one_particle_motion_z}
t\mapsto \pi_{(x,s)}(s+t)-x \qquad\text{ for }t\in[0,\infty)
\end{equation}
is a $D_{[0,\infty)}(\R)$-valued random variable,
with a law that does not depend on $z$.
We refer to this law as the one-particle motion of $\mc{A}$.

For a weave $\mc{A}$, let $\mc{F}_t$ be the $\sigma$-field generated by $\{\pi(s)\-\pi\in\mc{A},\,\sigma_\pi\leq  s\leq t\}$.
We say that $(\mc{F}_t)_{t\in\R}$ is the generated filtration of $\mc{A}$.

\begin{theorem}[Tightness criteria for homogeneous weaves]
\label{t:application_noncr}
Let $(\mc{A}_n)_{n\in\N}$ be a sequence of weaves,
with generated filtrations $(\mc{F}^n_t)_{t\in\R}$.
Assume the following:
\begin{enumerate}[series=theoremconditions]
\item 
For each $n\in\N$, $\mc{A}_n$ is homogeneous.
\item For each $n\in\N$, let $X_n$ be the one-particle motion of $\mc{A}_n$ from the origin. 
Suppose that $X_n$ takes values in $D_{[0,\infty)}(\R)$
and is strongly Markov with respect to the filtration $(\mc{F}^n_t)_{t\in[0,\infty)}$.
\item The laws of $(X_n)$ are tight on $D_{[0,\infty)}(\R)$, under J1 or M1.
\end{enumerate}
Then the laws of $(\mc{A}_n)$ are tight on $\mc{K}(\Pi^\uparrow)$ under M1 
and all weak limit points are non-crossing. 
\end{theorem}

We state also a more general, but more technical version of the above result.
We remark that it contains a slightly subtle point: 
the random path $\pi^n_z$, introduced within the statement of Theorem \ref{t:noncr_appliation_nonhom}, 
is not required to be an element of $\mc{A}_n$.

\begin{theorem}[Tightness criteria for inhomogeneous weaves]
\label{t:noncr_appliation_nonhom}
In the setting of Theorem \ref{t:application_noncr},
the same conclusion can be drawn if instead of conditions (i)-(iii) we have the following.
For each $n\in\N$ and for each deterministic $z\in\R^2$ there exists a random variable $\pi^n_z$ 
taking values in $\Pi^\uparrow$
(defined  on the same probability space as $\mc{A}_n$),
such that $\pi^n_z$ starts at $z$ and a.s.\ does not cross $\mc{A}_n$, and:
\begin{enumerate}[resume=theoremconditions]
\item
$\pi^n_z$ is almost surely continuous at $\sigma_{\pi^n_z}$ and is strongly Markov with respect to $(\mc{F}^n_t)$.
\item
For all $\eps>0$, $x\in\R$ and $T\in(0,\infty)$,
\begin{equation}
\label{eq:iii_iv_replacement}
\lim_{\de\to 0}
\sup_{n\in\N}\sup_{z\in\{x\}\times[-T,T]}
\P\l[\sup_{s\in[0,\de]}\l|\pi^n_z(\sigma_{\pi^n_z}+s)-\pi^n_z(\sigma_{\pi^n_z})\r|\geq\eps\r]=0.
\end{equation}
\end{enumerate}
\end{theorem}

\bpro
The proof of Theorems~\ref{t:application_noncr} and \ref{t:noncr_appliation_nonhom},
which we give together, 
will take up the remainder of Section~\ref{s:application_noncr}.
Let us first handle the connection between the two results, 
namely we show that (i)-(iii) implies both (iv) and (v).
To see this
take $\pi^n_z$ to be the one-particle motion of $\mc{A}_n$ from $z$,
which by (i) and our general comments above is well defined for all $z\in\R^2$,
with a distribution that does not depend on $z$.
Condition (iv) then follows from (ii) and translation invariance.
Condition (v) follows from (iii) and Lemma \ref{l:D_tightness}.
It therefore suffices to establish the conclusion of Theorem \ref{t:application_noncr}
under conditions (iv) and (v), which we assume hereon.

Let $T,\de,\eps\in(0,\infty)$ and $r\in\R$. 
For each $n\in\N$ we define a sequence of paths $\pi^n_m$
and stopping times $\tau^n_m$.
Set $\tau^n_0=-T$ and then for $m\geq 1$,
\begin{alignat*}{2}
\pi_m^n &= \pi_{(r+\eps,\tau_{m-1})}^n, \qquad
&&\tau_m^n = \inf\{t\geq \tau_{m-1}^n\-\pi_m^n(t)\leq r\text{ or }\pi_m^n(t)\geq r+2\eps\}.
\end{alignat*}
In words, the path $\pi^n_m$ begins at time $\tau^n_{m-1}$ at spatial location $r+\eps$.
The time $\tau^n_{m}$ occurs when $\pi^n_m$ exits the region $[r,r+2\eps]$,
upon which $\pi^n_{m+1}$ is born at space-time location $(r+\eps,\tau^n_m)$.
If $\tau^n_m$ is infinite, 
then the inductive definition terminates
and the sequences defined are taken to be finite.

\begin{figure}[t]
    \centering
    \includegraphics[width=0.5\textwidth]{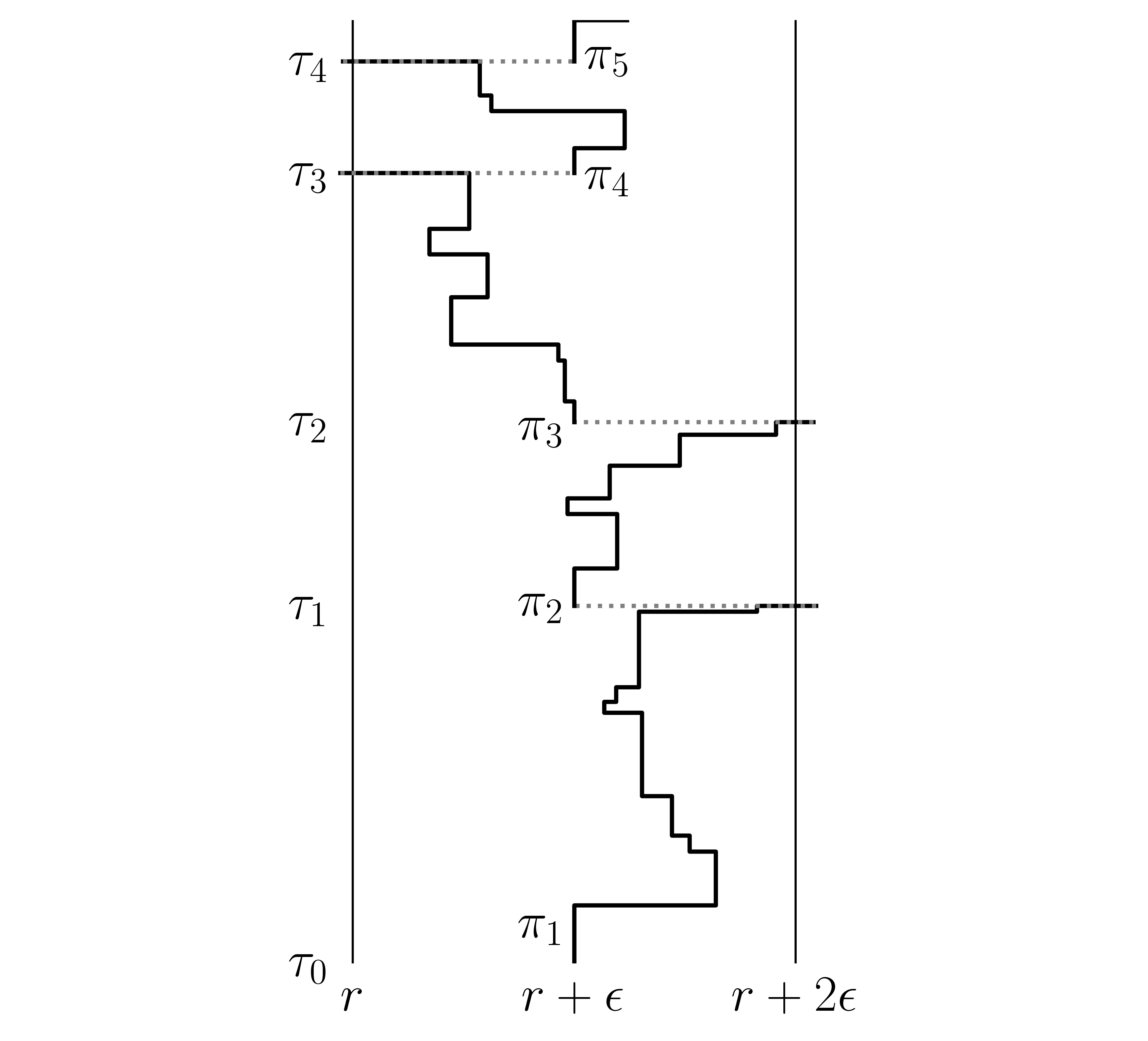} 
	\caption{A depiction of the sequences $(\tau^n_m)_{m\in\N}$ and $(\pi^n_m)_{m\in\N}$, with the superscript $n$ dropped. Time runs upwards and space is on the horizontal axis. Note that if a path $\pi$ is to cross $[r,r+2\eps]$ then it must pass through at least one of the grey dotted sections at time $\tau^n_m$. Moreover, any such crossing must be from left to right if the corresponding $\pi^n_m$ exits $[r,r+2\eps]$ via the right boundary, or must be from from right to left if the corresponding $\pi^n_m$ exits $[r,r+2\eps]$ via the left boundary. The event $E^n_\de$ ensures that the dotted sections are separated from one another in time by at least $\de$, which prevents the event $\mc{C}^{\rm M}_{T,\de,2\eps,0}\cap (\mc{A}_n\times\mc{A}_n)\neq\emptyset $ from occurring.}
	   \label{fig:applications_E_delta}
\end{figure}

Define also the events
\begin{equation*}
E_{\de,m}^n=\l\{\tau^n_{m}>\tau^n_{m-1}+\de\r\},\qquad
E_\de^n=\bigcap_{m=1}^{M_n} E_{\de,m}^n
\end{equation*}
where ${M_n}=\inf\{m\in\N\-\tau_m\geq T\}$.
We claim that, in the notation of \eqref{crotig},
\begin{equation}
\label{eq:ECWW}
\l\{\mc{C}^{\rm M}_{T,\de,2\eps,r}\cap (\mc{A}_n\times\mc{A}_n)\neq\emptyset\r\} \;\sw\; \Omega \sc E_\de^n.
\end{equation}
The reasoning behind equation \eqref{eq:ECWW} is explained in Figure \ref{fig:applications_E_delta},
which depicts the sequences $(\tau^n_m)$ and $(\pi^n_m)$.
In order to check \eqref{crotig} we must therefore check that the event $E^n_\de$
occurs with high probability, in a sense matching the limits appearing in \eqref{crotig}.

By \eqref{eq:iii_iv_replacement}
there exists $\de_0>0$ such that 
$$\inf_n\inf_{z\in\{r+\eps\}\times[-T,T]}\P\l[\sup_{s\in[0,\de_0]}\l|\pi^n_z(\sigma_{\pi^n_z}+s)-\pi^n_z(\sigma_{\pi^n_z})\r|<\eps\r]\geq\frac12,$$
which implies that, for all $m\in\N$,
$\inf_n\P\big[\tau^n_{m}-\tau^n_{m-1} \geq \de_0\big]\geq\frac12$.
The strong Markov property from condition (iv) gives that, for fixed $n$, the variables $ (\tau^n_{m}-\tau^n_{m-1})_{m\in\N}$
are independent.
Noting that $\tau^n_m-\tau^n_0=\sum_{i=1}^m \tau^n_m-\tau^n_{m-1}$,
it follows that $\tau^n_m$ is, for all $n$,
stochastically bounded below by $-T+\sum_1^m T_i$ where $\P[T_i\geq\de_0]\geq\frac12$. Hence $M_n$ is bounded above by a negative binomial distribution with parameters $N=\lceil 2T/\de_0 \rceil$ (the number of `successes') and $p=\frac12$. Noting that these parameters are independent of $n$, 
it follows that 
\begin{equation}
\label{eq:step_control}
\text{for all $\kappa>0$
there exists $K\in\N$ such that }
\sup_{n\in\N}\P[M_n\geq K]\leq\kappa.
\end{equation}
Let $\kappa>0$ and take $K\in\N$ as in \eqref{eq:step_control}, then
\begin{align}
\P[\Omega\sc E^n_\de]
&\leq \P[M_n\geq K] + \P\l[\bigcup_{m=1}^{M_n}\Omega\sc E^n_{\de,m}\text{ and }M_n<K\r] \notag \\
&\leq \kappa + \sum_{m=1}^{K}\P[\Omega\sc E^n_{\de,m}] \notag \\
&\leq \kappa + \sum_{m=1}^K \P\l[\sup_{s\in[0,\de]}\l|\pi^n_z(\sigma_{\pi^n_z}+s)-\pi^n_z(\sigma_{\pi^n_z})\r|\geq \eps\r]. \label{eq:En_control}
\end{align}
The first line of the above follows by elementary set algebra and conditioning. 
The second line follows from \eqref{eq:step_control} and the third from the definition of $E^n_{\de,m}$. 
Taking suprema over $n\in\N$ in \eqref{eq:En_control},
and then taking the $\limsup$ of both sides as $\de\to 0$, we obtain
\begin{equation}
\label{eq:En_control_2}
\limsup_{\de\to 0}\sup_{n\in\N}\P[\Omega\sc E^n_\de] \leq \kappa + \sum_{m=1}^K \limsup_{\de\to 0}\sup_{n\in\N} \P\l[\sup_{s\in[0,\de]}\l|\pi^n_z(\sigma_{\pi^n_z}+s)-\pi^n_z(\sigma_{\pi^n_z})\r|\geq \eps\r].
\end{equation}
We have from (v) that all terms within the summation on the right hand side are zero,
thus 
$\limsup_{\de\to 0}\sup_{n\in\N} \P[\Omega\sc E^n_\de]\leq \kappa.$ 
Using that $\kappa>0$ was arbitrary we obtain $\lim_{\de\to 0}\sup_{n\in\N} \P[\Omega\sc E^n_\de]=0$. In view of \eqref{eq:ECWW} and Theorem \ref{T:noncros}, this completes the proof.
\epro

\subsection{Heavy-tailed Poisson trees}
\label{s:application_poisson_trees}

In this section we use Theorem \ref{t:application_noncr}
to show tightness of a specific particle system, under its natural rescaling.
We treat a heavy-tailed version of the Poisson trees considered by \cite{FFW05, EFS17} and others,
in which
the motion of a single particle is within the natural domain of attraction of an $\alpha$-stable process, where $\alpha\in(0,2)$.
For brevity we show only tightness. We do not attempt to characterise the limit,
which will be a system of highly correlated non-crossing $\alpha$-stable processes.
In the case $\alpha=2$, 
after a diffusive rescaling the particles become coalescing Brownian motions that are independent before coalescence; 
under suitable conditions the limit in this case is known to be the Brownian web.
However in the $\alpha$-stable case one should expect that particles in the limit are dependent even before coalescence.

Fix $\alpha\in(0,2)$ and
let $\mu$ be a finite measure on $(0,\infty)$
such that 
\begin{equation}
\label{eq:stable_attractrion}
\lim_{R\to\infty} R^{\alpha}\int_R^\infty r\,\mu(dr)\;\in(0,\infty),
\end{equation}
for example $\mu(dr)=(1\wedge r^{-\alpha-2})\,dr$.
For each $n\in\N$
let $P_n$ denote a Poisson point process in $\R\times\R\times(0,\infty)$ with intensity measure
\begin{equation}
\label{eq:ppp}
n^{1/\alpha}\,dx \otimes n\,dt \otimes \mu_n(dr)
\end{equation}
where $\mu_n(A)=\mu\l(n^{1/\alpha} A\r)$.

Given a point $z=(x,t)\in\R^2$ we define a {\cadlag} path $\pi^n_{z}$ with time domain $[t,\infty)$
and initial point $(x,t)$, 
by specifying that $\pi^n_z$ remains constant except
at values of $t$ for which $\pi^n_z(t-)\in[x-r,x+r]$ for some $(x,t,r)\in P_n$;
at such times the path jumps and $\pi^n_z(t+)=x$.
We will show below that
$\pi_z^{n}$ is a compound Poisson process,
which by symmetry has zero drift.

\begin{figure}[t]
    \centering
    \includegraphics[width=0.5\textwidth]{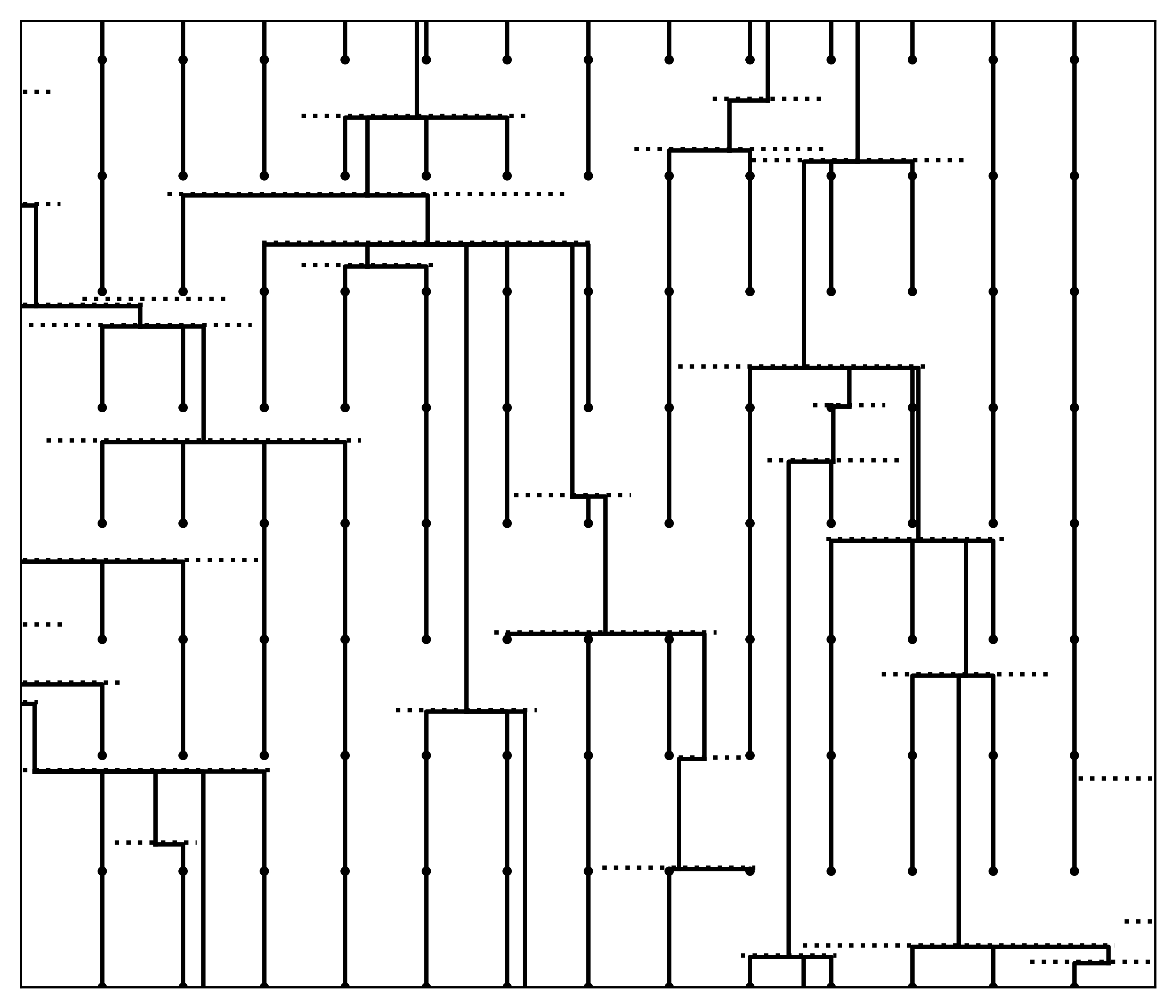} 
	\caption{Part of a Poisson tree, with one path beginning at each black dot and a square border encasing the space-time window displayed. Time runs upwards and space is on the horizontal axis. Events $(x,t,r)\in P_n$ are shown as dotted sections of the form $[x-r,x+r]\times\{t\}$. The realization of the underlying Poisson process $P_n$ has been chosen to make a clear picture, rather than to demonstrate heavy-tailed behaviour.}
	   \label{fig:applications_poisson_tree}
\end{figure}

We define
$$W_n=\{\pi^n_z\-z\in\R^2\}$$
and let $\mc{W}_n$ be the closure of $W_n$ in $\Pi^\uparrow$.
See Figure \ref{fig:applications_poisson_tree} for a graphical depiction of $W_n$ and $P_n$.
Equation \eqref{eq:ppp} corresponds to a space-time rescaling in which, at stage $n$, 
we speed up time by a factor $n$ and compress space by a factor $n^{1/\alpha}$.
In particular, 
for fixed but arbitrary $z=(x,t)\in\R^2$ and $n\in\N$, 
for $s\in(0,\infty)$ the processes
\bc
\label{eq:poisson_trees_scaling_paths}
&s\mapsto \pi^{n}_{z}(t+s)-x \quad \text{ and }\quad s\mapsto n^{-1/\alpha}\l(\pi^{1}_{z}(t+ns)-x\r)
\quad \text{ have the same law.} 
\ec
Moreover, 
since \eqref{eq:ppp} is invariant in law under deterministic translations of space-time,
the law in \eqref{eq:poisson_trees_scaling_paths} does not depend on $z$.

The remainder of the present section will apply Theorem \ref{t:application_noncr} to the sequence $(\mc{W}_n)$.
We must first show that $\mc{W}_n$ is a weave.
Then we must check conditions (i)-(iii) of Theorem \ref{t:application_noncr};
conditions (i) and (ii) are essentially immediate and condition (iii) will come from showing 
that the limiting one-particle motion is $\alpha$-stable.

\begin{remark}
For arbitrary $Z_n\sw\R^2$ 
we have that $A_n=\{\pi^n_z\-z\in Z_n\}\sw W_n\sw\mc{W}_n$.
Setting $\mc{A}_n$ to be the closure of $A_n$,
tightness of the sequence $(\mc{A}_n)$
follows immediately from that of $(\mc{W}_n)$.
For example, we might use paths begun on rescalings of the square lattice $\Z^2_n=\{(x,t)\in\R^2\-n^{1/\alpha} x\in\Z,\;nt\in\Z\}$
as depicted in Figure \ref{fig:applications_poisson_tree}.
\end{remark}

We begin with the one-particle motion.
More precisely, we will show that the 
process $X_n(s)=\pi^{(n)}_{(0,0)}(s)$ is a compound Poisson process.
In view of \eqref{eq:poisson_trees_scaling_paths} it suffices to consider $n=1$.
Let us note at the outset that condition \eqref{eq:stable_attractrion}
combined with the requirement that $\mu$ be a finite measure
implies that $\int_0^\infty r\,\mu(dr)<\infty$.

We refer to each $(x,t,r)\in P_n$ as an event and to the space-time region $[x-r,x+r]\times\{t\}$ as 
being affected by the event. We similarly say that when $X_n(t-)\in[x-r,x+r]$,
the path $X_n$ is affected by the event $(x,t,r)$.
Taking $n=1$,
if the current location of $X_1=\pi^{(1)}_{(0,0)}$ is $y\in\R$ then it is affected by events at rate
\begin{align}
\label{eq:event_hit_rate}
\int_0^\infty \int_\R \1_{\{y\in[x-r,x+r]\}}\,dx\,\mu(dr) 
=2 \int_0^{\infty} r\,\mu(dr)
\end{align}
which is finite.
Upon being affected by an event, 
the resulting jump changes the spatial location of $X_n$ by (addition of) a random variable $J$ for which
\begin{align}
\P[J\geq R] &=
\frac{1}{K}  \int_0^\infty \int_\R \1_{\{y\in[x-r,x-R]\}}\,dx\,\mu(dr) \notag \\
&= \frac{1}{K} \int_0^\infty (r-R)\vee 0 \;\mu(dr) \notag \\
&= \frac{1}{K} \int_R^\infty r \,\mu(dr),\label{eq:event_hit_displacement}
\end{align}
where $R\in(0,\infty)$ and $K=2\int_0^\infty r\,\mu(dr)$.
Note that by symmetry for $R\geq 0$ we have  $\P[J\geq R]=\P[J\leq -R]$,
so \eqref{eq:event_hit_displacement} characterises the distribution of $J$.

\begin{prop}[Tightness of heavy-tailed Poisson trees]
\label{p:poisson_trees_tight}
The laws of $(\mc{W}_n)$ are tight on $\mc{K}(\Pi^{\uparrow})$ under M1 
and all weak limit points are non-crossing.
\end{prop}
\bpro
We first justify that $\mc{W}_n$ is a weave, for each $n\in\N$.
With $n\in\N$ fixed, it is clear that $\mc{W}_n$ is a subset of $\Pi^\uparrow$
and that each $z\in\R^2$ satisfies $z\in\ov{\pi}$ for some $\pi\in\mc{W}_n$.
To see that $\mc{W}_n$ is a weave, it remains only to show that $\mc{W}_n$ is almost surely compact and non-crossing, for each $n\in\N$.

A similar calculation to \eqref{eq:event_hit_rate}
shows that for each $L,T\in(0,\infty)$ there exists, almost surely, 
a random $\de>0$ (depending on $L$, $T$ and $n$) such that the regions $[x-r,x+r]\times[t-\de,t+\de]$
are disjoint for all $(x,t,r)\in P_n$ with $(x,t)\in[-L,L]\times[-T,T]$.
Note that paths in $W_n$ remain constant (in space) outside of such regions.
With this in hand it is easily seen via Theorem \ref{T:ArzAsc} that $\mc{W}_n=\ov{W_n}$ is almost surely compact, 
under both J1 and M1.
Moreover, noting that $W_n$ is non-crossing,
Proposition~\ref{P:collis} gives that if two paths $\pi,\pi'\in\mc{W}_n$ were to cross, then they must collide
at time $\sigma_\pi\vee\sigma_{\pi'}$.
However, this eventuality cannot occur due to $\de>0$.
Hence $\mc{W}_n$ is almost surely non-crossing.

To complete the proof we will apply Theorem \ref{t:application_noncr} to the sequence $(\mc{W}_n)_{n\in\N}$,
for which we must check conditions (i)-(iii) of that theorem.
Conditions (i) and (ii) follow immediately from standard independence properties and translation-invariance of Poisson point processes, so it remains only to check condition (iii).

We have seen above that the one-particle motion of $\mc{W}_n$
from the origin is a compound Poisson process with jump rate \eqref{eq:event_hit_rate} 
and jump distribution $J$ characterised by \eqref{eq:event_hit_displacement}.
Theorem 4.5.2 of \cite{Whi02} gives that
a \emph{discrete} time random walk with i.i.d.~increments having distribution $J$ is within the normal domain of attraction of an $\alpha$-stable process, with time rescaled by a factor $n$ and space rescaled by $n^{-1/\alpha}$. 
Theorem 4.5.3 of \cite{Whi02}
gives that such a walk converges in law, as a {\cadlag} process under the J1 (and hence also M1) topology, to an $\alpha$-stable process.
For brevity we do not calculate the scale parameter of the limiting $\alpha$-stable process here
but it may be found via formulae therein; by symmetry the skewness and shift parameters are both zero.
Noting \eqref{eq:event_hit_rate},
we are in fact concerned with a random walk in continuous time, that is, a compound Poisson process.
The two cases differ only by a strictly increasing piecewise linear time change,
a time change which in the limit converges almost surely, in the locally uniform sense, to a linear time change.
Thus, by Theorem 13.2.2 of \cite{Whi02},
the same results apply to the sequence $n\mapsto X_n(\cdot)$.
This establishes condition (iii) of Theorem \ref{t:application_noncr},
which completes the proof.
\epro

\section{Proofs of the tightness criteria}\label{S:tightproof}

\subsection{The Hausdorff topology}

Let $\Xc$ be a metrisable space, let $d$ be any metric generating the topology, and let $\Ki_+(\Xc)$ be the space of nonempty compact subsets of $\Xc$, equipped with the Hausdorff topology. The following lemma, which we cite from \cite[Lemma~B.1]{SSS14}, provides a link between convergence of compact sets in the topology on $\Ki_+(\Xc)$, and convergence of their elements in the topology on the underlying space $\Xc$. This will be used many times in what follows. Note that it suffices to prove the first and last inclusion in (\ref{Haulim}) to conclude that all sets therein are equal.


\bl[Convergence in the Hausdorff topology]
Let\label{L:Hauconv} $K_n,K\in\Ki_+(\Xc)$. Then $K_n\to K$ in the Hausdorff topology if and only if there exists a $C\in\Ki_+(\Xc)$ such that $K_n\sub C$ for all $n$ and
\bc\label{Haulim}
K&\sub&\dis\big\{x\in\Xc:\exists x_n\in K_n\mbox{ s.t.\ }x_n\to x\big\}\\[5pt]
&\sub&\dis\big\{x\in\Xc:\exists x_n\in K_n
\mbox{ s.t.\ $x$ is a cluster point of } (x_n)_{n\in\N}\big\}\sub K.
\ec
\el

As a direct consequence of Lemma~\ref{L:Hauconv} one obtains:

\bl[Convergence of continuous images]
Let\label{L:contim} $\Xc,\Yi$ be metrisable and let $\psi\cn\Xc\to\Yi$ be continuous. Then $K_n\to K$ in $\Ki_+(\Xc)$ implies $\psi(K_n)\to\psi(K)$ in $\Ki_+(\Yi)$.
\el

We will moreover need the following lemma, which we cite from \cite[Lemma~2.12]{FS25}. Note that the lemma implies that if $\Xc$ is compact, then so is $\Ki_+(\Xc)$.

\bl[Compactness in the Hausdorff topology]
A\label{L:Haucomp} set $\Ai\sub\Ki_+(\Xc)$ is precompact if and only if there exists a $C\in\Ki_+(\Xc)$ such that $K\sub C$ for all $K\in\Ai$.
\el

If $\pi\in\Pi$ is a path and $\ov\pi$ is its filled graph, then we set
\be\label{pinc}
\pi^{\li 2\re}:=\big\{(z,z')\in\pi\times\pi:z\pre z'\big\}
\quand
\ov\pi^{\li 2\re}:=\big\{(z,z')\in\ov\pi\times\ov\pi:z\pre z'\big\},
\ee
where $\pre$ denotes the total order on $\pi$ or $\ov\pi$, respectively. It is not hard to see that $\pi^{\li 2\re}$ and $\ov\pi^{\li 2\re}$ are compact subsets of $\R^2_{\rm c}\times\R^2_{\rm c}$, equipped with the product topology. The following lemma is a consequence of \cite[Thm~1.9]{FS25}.

\bl[Convergence in the J1 and M1 topologies]
Let\label{L:JMconv} $\pi_n,\pi\in\Pi$. Then one has $\pi_n\to\pi$ in the J1 topology if and only if $\pi_n^{\li 2\re}\to\pi^{\li 2\re}$ in the Hausdorff topology on $\Ki_+(\R^2_{\rm c}\times\R^2_{\rm c})$. Similarly, $\pi_n\to\pi$ in the M1 topology if and only if $\ov\pi_n^{\li 2\re}\to\ov\pi^{\li 2\re}$ in the Hausdorff topology.
\el

The already mentioned \cite[Thm~1.9]{FS25} also yields the following statement, which can alternatively be derived from Lemma~\ref{L:JMconv} using Lemma~\ref{L:contim}.

\bl[Hausdorff convergence of paths]
If\label{L:Haupi} $\pi_n\to\pi$ in the J1 topology on $\Pi$, then $\pi_n\to\pi$ in the Hausdorff topology on $\Ki_+(\R^2_{\rm c})$. Similarly, if $\pi_n\to\pi$ in the M1 topology on $\Pi$, then $\ov\pi_n\to\ov\pi$ in the Hausdorff topology on $\Ki_+(\R^2_{\rm c})$.
\el

\subsection{Compactness criteria}

Let $\dr$ be any metric generating the topology on the extended real line $\ov\R$. Let $\dr(x,A):=\inf_{y\in A}d(x,y)$ denote the distance of a point $x$ to a subset $A\sub\ov\R$. For any path $\pi\in\Pi$, we define moduli of continuity by
\bc
\dis m_{T,\de}(\pi)&:=&\dis\sup\big\{\dr(x,y):(x,y)\in\De^2_{T,\de}(\pi)\big\},\\[5pt]
\dis m^{\rm J}_{T,\de}(\pi)&:=&\dis\sup\big\{\dr(y,\{x,z\}):(x,y,z)\in\De^3_{T,\de}(\pi)\big\},\\[5pt]
\dis m^{\rm M}_{T,\de}(\pi)&:=&\dis\sup\big\{\dr(y,[x,z]):(x,y,z)\in\De^3_{T,\de}(\pi)\big\}.
\ec
Then \cite[Thms 1.11 and 1.12]{FS25} tell us the following.

\bt[Compactness criteria]
A\label{T:ArzAsc} set $\Ai\sub\Pi_{\rm c}$ is precompact if and only if
\be\label{Arz1}
\lim_{\de\to 0}\sup_{\pi\in\Ai}m_{T,\de}(\pi)=0\quad\forall T>0.
\ee
A set $\Ai\sub\Pi$ is precompact with respect to the J1 topology if and only if
\be\label{Arz2}
\lim_{\de\to 0}\sup_{\pi\in\Ai}m^{\rm J}_{T,\de}(\pi)=0\quad\forall T>0.
\ee
A set $\Ai\sub\Pi$ is precompact with respect to the M1 topology if and only if
\be\label{Arz3}
\lim_{\de\to 0}\sup_{\pi\in\Ai}m^{\rm M}_{T,\de}(\pi)=0\quad\forall T>0.
\ee
\et

Theorem~\ref{T:ArzAsc} and Lemma~\ref{L:Haucomp} allow us to characterise the precompact subsets of $\Ki_+(\Pi_{\rm c})$, respectively of $\Ki_+(\Pi)$, with respect to the J1 and M1 topologies. We can use this to prove some first tightness criteria for the laws of random compact sets of paths, from which we will ultimately derive Theorems \ref{T:ctight} and  \ref{T:tight}. We start with continuous paths.

\bt[Tightness criterion for sets of continuous paths]
Let\label{T:Rctight} $(\Ai_\ga)_{\ga\in\Ga}$ be a family of random variables with values in $\Ki_+(\Pi_{\rm c})$. Then the laws $(\mu_\ga)_{\ga\in\Ga}$ with $\mu_\ga:=\P[\Ai_\ga\in\,\cdot\,]$ are tight with respect to the topology on $\Pi_{\rm c}$ if and only if
\be\label{Rctight}
\lim_{\de\to 0}\sup_{\ga\in\Ga}\P\big[\sup_{\pi\in\Ai_\ga}m_{T,\de}(\pi)>\eps\big]=0\quad\forall T,\eps>0.
\ee
\et

\bpro
Assume that the laws $(\mu_\ga)_{\ga\in\Ga}$ are tight. Then for each $\eta>0$, there exists a compact $\Cb_\eta\sub\Ki_+(\Pi_{\rm c})$ such that
\be
\sup_{\ga\in\Ga}\P\big[\Ai_\ga\in\Cb_\eta\big]\geq 1-\eta.
\ee
By Lemma~\ref{L:Haucomp}, there exists a compact $\Ci_\eta\sub\Pi_{\rm c}$ such that $\Ai\in\Cb_\eta$ implies $\Ai\sub\Ci_\eta$. Thus, for each $\eta>0$, there exists a compact $\Ci_\eta\sub\Pi_{\rm c}$ such that
\be\label{AC}
\P\big[\Ai_\ga\sub\Ci_\eta\big]\geq\P\big[\Ai_\ga\in\Cb_\eta\big]\geq 1-\eta\qquad(\ga\in\Ga).
\ee
By (\ref{Arz1}) of Theorem~\ref{T:ArzAsc},
\be\label{w0}
\lim_{\de\to 0}w_{T,\de}(\eta)=0\quad(T>0)
\quad\mbox{with}\quad
w_{T,\de}(\eta):=\sup_{\pi\in\Ci_\eta}m_{T,\de}(\pi)\quad(T,\de>0).
\ee
Fix $T,\eps>0$. It follows that for each $\eta>0$, there exists a $\de>0$ such that $w_{T,\de}(\eta)\leq\eps$ and hence, by (\ref{AC}) and the definition of $w_{T,\de}(\eta)$,
\be
\sup_{\ga\in\Ga}\P\big[\sup_{\pi\in\Ai_\ga}m_{T,\de}(\pi)>\eps\big]
\leq\sup_{\ga\in\Ga}\P\big[\Ai_\ga\not\sub\Ci_\eta\big]\leq\eta,
\ee
proving (\ref{Rctight}).

Assume, conversely, that (\ref{Rctight}) holds. Formula (\ref{Rctight}) implies that for each $\eta>0$ and integers $T,n\geq 1$, there exists a $\de_n(T,\eta)>0$ such that
\be
\sup_{\ga\in\Ga}\P\big[\sup_{\pi\in\Ai_\ga}m_{T,\de_n(T,\eta)}(\pi)>n^{-1}\big]\leq\eta 2^{-T-n}.
\ee
Summing over $T$ and $n$, it follows that
\be
\P\big[\exists T,n\geq 1\mbox{ s.t.\ }\sup_{\pi\in\Ai_\ga}m_{T,\de_n(T,\eta)}(\pi)>n^{-1}\big]\leq\eta\qquad(\ga\in\Ga).
\ee
Setting
\be
\Ci_\eta:=\big\{\pi\in\Pi_{\rm c}:m_{T,\de_n(T,\eta)}(\pi)\leq n^{-1}\ \forall T,n\geq 1\big\},
\ee
our previous formula says that
\be
\P\big[\Ai_\ga\sub\Ci_\eta\big]\geq 1-\eta\qquad(\ga\in\Ga).
\ee
We observe that $\de\leq\de_n(T,\eta)$ implies $m_{T,\de}(\pi)\leq m_{T,\de_n(T,\eta)}(\pi)$ and hence by the definition of $\Ci_\eta$
\be
\sup_{\pi\in\Ci_\eta}m_{T,\de}(\pi)\leq n^{-1}\qquad\forall T\geq 1,\ \de\leq\de_n(T,\eta).
\ee
It follows that
\be
\limsup_{\de\to 0}\sup_{\pi\in\Ci_\eta}m_{T,\de}(\pi)\leq n^{-1}\quad\forall T,n\geq 1,
\ee
which by (\ref{Arz1}) of Theorem~\ref{T:ArzAsc} implies that $\Ci_\eta$ is precompact. Letting $\ov\Ci_\eta$ denote its closure, we have for each $\eta>0$ found a compact set $\ov\Ci_\eta\sub\Pi_{\rm c}$ such that
\be
\sup_{\ga\in\Ga}\P\big[\Ai_\ga\sub\ov\Ci_\eta\big]\geq 1-\eta.
\ee
By Lemma~\ref{L:Haucomp}, this implies that the laws $(\mu_\ga)_{\ga\in\Ga}$ are tight.
\epro

In analogy with Theorem~\ref{T:Rctight} we obtain the following tightness criteria for compact sets of {\cadlag} paths.

\bt[Tightness criteria for sets of {\cadlag} paths]\hspace{5pt}
Let\label{T:Rtight} $(\Ai_\ga)_{\ga\in\Ga}$ be a family of random variables with values in $\Ki_+(\Pi)$. Then the laws $(\mu_\ga)_{\ga\in\Ga}$ with $\mu_\ga:=\P[\Ai_\ga\in\,\cdot\,]$ are tight with respect to the J1 topology on $\Pi$ if and only if
\be\label{Rtight}
\lim_{\de\to 0}\sup_{\ga\in\Ga}\P\big[\sup_{\pi\in\Ai_\ga}m^{\rm J}_{T,\de}(\pi)\geq\eps\big]=0\quad\forall T,\eps>0.
\ee
An analogue statement holds for the M1 topology, with $m^{\rm J}_{T,\de}$ replaced by $m^{\rm M}_{T,\de}$.
\et

\bpro
The proof is completely the same as the proof of Theorem~\ref{T:Rctight}, using (\ref{Arz2}) and (\ref{Arz3}) of Theorem~\ref{T:ArzAsc} instead of (\ref{Arz1}).
\epro

\subsection{Tightness criteria}

In this subsection we derive Theorems \ref{T:ctight} and \ref{T:tight} from Theorems \ref{T:Rctight} and \ref{T:Rtight}. Fix $T,\de>0$. Recall the definitions of the sets $\De^2_{T,\de}(\pi)$ and $\De^3_{T,\de}(\pi)$ in (\ref{De2}) and (\ref{De3}). For each $\pi\in\Pi$, we set
\bc
\dis\De^+_{T,\de}(\pi)&:=&\dis\big\{(x,y)\in\De^2_{T,\de}(\pi):x<y\big\},\\[5pt]
\dis\De^-_{T,\de}(\pi)&:=&\dis\big\{(x,y)\in\De^2_{T,\de}(\pi):x>y\big\},\\[5pt]
\dis\De^{++}_{T,\de}(\pi)&:=&\dis\big\{(x,y,z)\in\De^3_{T,\de}(\pi):x<y<z\big\},\\[5pt]
\dis\De^{+-}_{T,\de}(\pi)&:=&\dis\big\{(x,y,z)\in\De^3_{T,\de}(\pi):x<y>z\big\},\\[5pt]
\dis\De^{-+}_{T,\de}(\pi)&:=&\dis\big\{(x,y,z)\in\De^3_{T,\de}(\pi):x>y<z\big\},\\[5pt]
\dis\De^{--}_{T,\de}(\pi)&:=&\dis\big\{(x,y,z)\in\De^3_{T,\de}(\pi):x>y>z\big\}.
\ec
For each $\eta>0$ and $\star\in\{2,+,-\}$, we define
\be
\Ti^\star_{T,\de,\eta}:=\big\{\pi\in\Pi:\exists(x,y)\in\De^\star_{T,\de}\mbox{ s.t.\ }\dr(x,y)>\eta\big\},
\ee
and for $\star\in\{3,++,+-,-+,--\}$, we define
\be
\Ti^\star_{T,\de,\eta}:=\big\{\pi\in\Pi:\exists(x,y,z)\in\De^\star_{T,\de}\mbox{ s.t.\ }\dr(x,y),\dr(y,z)>\eta\big\}.
\ee
Finally, we set
\be
\Ti^{\rm J}_{T,\de,\eta}:=\Ti^{++}_{T,\de,\eta}\cup\Ti^{--}_{T,\de,\eta}
\quand
\Ti^{\rm M}_{T,\de,\eta}:=\Ti^{+-}_{T,\de,\eta}\cup\Ti^{-+}_{T,\de,\eta}.
\ee

\bl[Alternative tightness criteria]
Condition\label{L:alt} (\ref{Rctight}) is equivalent to
\be\label{cS}
\lim_{\de\to 0}\sup_{\ga\in\Ga}\P\big[\Ai_\ga\cap\Ti^2_{T,\de,\eta}\neq\emptyset\big]=0\quad\forall T,\eta>0.
\ee
Condition (\ref{Rtight}) is equivalent to
\be\label{JS}
\lim_{\de\to 0}\sup_{\ga\in\Ga}\P\big[\Ai_\ga\cap(\Ti^{\rm J}_{T,\de,\eta}\cup\Ti^{\rm M}_{T,\de,\eta})\neq\emptyset\big]=0\quad\forall T,\eta>0.
\ee
Condition (\ref{Rtight}) with $m^{\rm J}_{T,\de}$ replaced by $m^{\rm M}_{T,\de}$ is equivalent to
\be\label{MS}
\lim_{\de\to 0}\sup_{\ga\in\Ga}\P\big[\Ai_\ga\cap\Ti^{\rm M}_{T,\de,\eta}\neq\emptyset\big]=0\quad\forall T,\eta>0.
\ee
\el

\bpro
We have $\sup_{\pi\in\Ai_\ga}m_{T,\de}(\pi)>\eta$ if and only if there exist $\pi\in\Ai_\ga$ and $(x,y)\in\De^2_{T,\de}(\pi)$ such that $\dr(x,y)>\eta$, so (\ref{Rctight}) is clearly equivalent to (\ref{cS}).

We have $\sup_{\pi\in\Ai_\ga}m^{\rm J}_{T,\de}(\pi)>\eta$ if and only if there exist $\pi\in\Ai_\ga$ and $(x,y,z)\in\De^3_{T,\de}(\pi)$ such that $\dr(x,y)\wedge\dr(y,z)>\eta$, so (\ref{Rtight}) is clearly equivalent to (\ref{JS}).

We have $\sup_{\pi\in\Ai_\ga}m^{\rm M}_{T,\de}(\pi)>\eta$ if and only if there exist $\pi\in\Ai_\ga$ and $(x,y,z)\in\De^3_{T,\de}(\pi)$ such that $\dr(y,[x,z])>\eta$. Here $\dr(y,[x,z])>\eta$ is equivalent to
\be
\dr(x,y)\wedge\dr(y,z)>\eta\quad\mbox{and either }x,z<y\mbox{ or }x,z>y,
\ee
so (\ref{Rtight}) with $m^{\rm J}_{T,\de}$ replaced by $m^{\rm M}_{T,\de}$ is equivalent to (\ref{MS}).
\epro

In what follows, it will be convenient to make a concrete choice for the metric $\dr$ on $\ov\R$. We choose
\be
\dr(x,y):=\big|\phi(x)-\phi(y)\big|\qquad(x,y\in\ov\R)
\quad\mbox{with}\quad
\phi(x):=\frac{x}{\sqrt{1+x^2}}\quad(x\in\R),
\ee
and $\phi(\pm\infty):=\pm 1$. We observe that
\be\label{Lip}
\dr(x,y)=\int_x^y(1+z^2)^{-3/2}\,\di z<y-x\qquad(x,y\in\ov\R,\ x<y).
\ee
For $\eps>0$, 
we choose $k_\pm(\eps)\in\Z$ with $k_-(\eps)<0<k_+(\eps)$ such that
\be
\dr\big(\pm\infty,k_\pm(\eps)\eps\big)<\eps.
\ee
We need the following simple lemma.

\bl[Jumping over intervals]
For\label{L:over} each $\eps>0$ and $x,y\in\ov\R$ such that $x<y$ and $\dr(x,y)>2\eps$, there exists a $k\in\Z$ with $k_-(\eps)\leq k\leq k_+(\eps)-1$ such that $x\leq k\eps$ and $y\geq(k+1)\eps$.
For each $\eps>0$ and $x,y,z\in\ov\R$ such that $x<y<z$ and $\dr(x,y),\dr(y,z)>2\eps$, there exists a $k\in\Z$ with $k_-(\eps)\leq k\leq k_+(\eps)-3$ such that $x\leq k\eps$, $(k+1)\eps\leq y\leq(k+2)\eps$, and $(k+3)\eps\leq z$.
\el

\bpro
Fix $\eps>0$. If $x,y\in\ov\R$ satisfy $x<y$ and $\dr(x,y)>2\eps$, then by (\ref{Lip}) $y-x>2\eps$ and hence the set
\be
K:=\big\{k\in\Z:x\leq k\eps\mbox{ and }y\geq(k+1)\eps\big\}
\ee
is nonempty. We claim that if $k\in K$ satisfies $k<k_-(\eps)$, then also $k_-(\eps)\in K$. Indeed $k\in K$ and $k<k_-(\eps)$ imply $x\leq k\eps<k_-(\eps)\eps$ while the observation that
\be
\dr(-\infty,y)\geq\dr(x,y)>2\eps>\dr\big(-\infty,k_-(\eps)\eps\big)+\eps
>\dr\big(-\infty,(k_-(\eps)+1)\eps\big)
\ee
implies $y\geq(k_-(\eps)+1)\eps$, completing the proof that $k_-(\eps)\in K$. By the same argument, if $k\in K$ satisfies $k>k_+(\eps)-1$, then also $k_+(\eps)-1\in K$ and the first claim of the lemma follows.

The proof of the second claim goes a bit differently. Assume that $x,y\in\ov\R$ satisfy $x<y$ and $\dr(x,y)>2\eps$. Choose $k\in\Z$ such that $(k+1)\eps\leq y\leq(k+2)\eps$. Since $\dr(x,y)>2\eps$, (\ref{Lip}) tells us that $y-x>2\eps$ and hence $x\leq k\eps$. The same argument gives $(k+3)\eps\leq z$. Since $y-(k+1)\eps\leq\eps$ (\ref{Lip}) tells us that $\dr((k+1)\eps,y)<\eps$. Combining this with the inequality $\dr(-\infty,y)\geq\dr(x,y)>2\eps$ we see that $\dr(-\infty,(k+1)\eps)>\eps$ and hence $k_-(\eps)<k+1$. The same argument gives $k+2<k_+(\eps)$.
\epro

We next compare sets of the form $\Ti^\star_{T,\de,\eta}$ with the sets $\Si^\star_{T,\de,\eta,r}$ defined in (\ref{Sipm}).

\bl[Comparison of sets of paths]
For\label{L:setcom} $\star\in\{+,-,+-,-+\}$ one has
\be\label{setcom1}
\Si^\star_{T,\de,\eps,r}\sub\Ti^\star_{T,\de,\eta}\quad\mbox{with}\quad\eta:=\dr(r,r+\eps)\qquad(T,\de,\eps>0,\ r\in\R)
\ee
and
\be\label{setcom2}
\Ti^\star_{T,\de,2\eps}\sub\bigcup_{k=k_-(\eps)}^{k_+(\eps)-1}\Si^\star_{T,\de,\eps,k\eps}\qquad(T,\de,\eps>0).
\ee
For $\star\in\{++,--\}$ one has
\be\label{setcom3}
\Si^\star_{T,\de,\eps,r}\sub\Ti^\star_{T,\de,\eta}\quad\mbox{with}\quad\eta:=\dr(r,r+\eps)\wedge\dr(r+2\eps,r+3\eps)\qquad(T,\de,\eps>0,\ r\in\R)
\ee
and
\be\label{setcom4}
\Ti^\star_{T,\de,2\eps}\sub\bigcup_{k=k_-(\eps)}^{k_+(\eps)-3}\Si^\star_{T,\de,\eps,k\eps}\qquad(T,\de,\eps>0).
\ee
\el

\bpro
We first prove (\ref{setcom1}). If $\pi\in\Si^+_{T,\de,\eps,r}$, then there exists a pair $(x,y)\in\De^2_{T,\de}(\pi)$ with $x<y$, $x\leq r$, and $r+\eps\leq y$. Then $\dr(x,y)\geq\dr(r,r+\eps)$ and hence $\pi\in\Ti^+_{T,\de,\eta}$ with $\eta:=\dr(r,r+\eps)$. This proves (\ref{setcom1}) for $\star=+$. The same argument with the roles of $x$ and $y$ reversed yields (\ref{setcom1}) for $\star=-$. The arguments for $\star=+-$ and $\star=-+$ are also very similar, except that there are now three points $x,y,z$ of which $x,z$ lie on one side of the interval $[r,r+\eps]$ while $y$ lies on the other side.

The proof of (\ref{setcom3}) is also similar. If $\pi\in\Si^+_{T,\de,\eps,r}$, then there exists a triple $(x,y,z)\in\De^3_{T,\de}(\pi)$ with $x<y<z$, $x\leq r$, $r+\eps\leq y\leq r+2\eps$, and $r+3\eps\leq z$. Defining $\eta$ as in (\ref{setcom3}) we then have $\dr(x,y)\geq\dr(r,r+\eps)\geq\eta$ and $\dr(y,z)\geq\dr(r+2\eps,r+3\eps)\geq\eta$ which implies that $\Ti^\star_{T,\de,\eta}$. This proves (\ref{setcom3}) for $\star=++$. The proof for $\star=--$ is the same with the roles of $x$ and $z$ reversed.

We next prove (\ref{setcom2}). If $\pi\in\Ti^+_{T,\de,2\eps}$ then there exist $(x,y)\in\De^2_{T,\de}(\pi)$ with $x<y$ and $\dr(x,y)>2\eps$. By Lemma~\ref{L:over} there then exists a $k\in\Z$ with $k_-(\eps)\leq k<k_+(\eps)$ such that $x\leq k\eps$ and $y\geq(k+1)\eps$. Then $\pi\in\Si^+_{T,\de,\eps,k\eps}$. This proves (\ref{setcom2}) for $\star=+$. The same argument with the roles of $x$ and $y$ reversed yields (\ref{setcom2}) for $\star=-$.
For $\star=+-$ we argue as follows. If $\pi\in\Ti^{+-}_{T,\de,2\eps}$ then there exist $(x,y,z)\in\De^3_{T,\de}(\pi)$ with $x,z<y$ and $\dr(x,y),\dr(z,y)>2\eps$. Then $x\vee z<y$ and $\dr(x\vee z,y)>2\eps$, so by Lemma~\ref{L:over} there then exists a $k\in\Z$ with $k_-(\eps)\leq k<k_+(\eps)$ such that $x\vee z\leq k\eps$ and $y\geq(k+1)\eps$, proving that $\pi\in\Si^{+-}_{T,\de,\eps,k\eps}$. The proof of (\ref{setcom2}) for $\star=-+$ is the same, using $x\wedge z$ instead of $x\vee z$.

It remains to prove (\ref{setcom4}). If $\pi\in\Ti^{++}_{T,\de,2\eps}$ then there exist $(x,y,z)\in\De^3_{T,\de}(\pi)$ with $x<y<z$ and $\dr(x,y),\dr(z,y)>2\eps$. By Lemma~\ref{L:over} there then exists a $k\in\Z$ with $k_-(\eps)\leq k\leq k_+(\eps)-3$ such that $x\leq k\eps$, $(k+1)\eps\leq y\leq(k+2)\eps$, and $(k+3)\eps\leq z$. Then $\pi\in\Si^{++}_{T,\de,\eps,k\eps}$. This proves (\ref{setcom4}) for $\star=++$. The argument for $\star=--$ is the same with the roles of $x$ and $z$ reversed.
\epro

\bpro[of Theorem~\ref{T:ctight}]
By Theorem~\ref{T:Rctight} and Lemma~\ref{L:alt} it suffices to prove that (\ref{ctight}) is equivalent to (\ref{cS}). By formula (\ref{setcom1}) of Lemma~\ref{L:setcom}, for all $T,\de,\eps>0$ and $r\in\R$,
\be
\P\big[\Si^2_{T,\de,\eps,r}\cap\Ai_\ga\neq\emptyset\big]
\leq\P\big[\Ti^2_{T,\de,\eta}\cap\Ai_\ga\neq\emptyset\big]
\quad\mbox{with}\quad\eta:=\dr(r,r+\eps)\qquad(T,\de,\eps>0,\ r\in\R),
\ee
from which we see that (\ref{cS}) implies (\ref{ctight}). By formula (\ref{setcom2}) of Lemma~\ref{L:setcom}, for all $T,\de,\eps>0$,
\be
\P\big[\Ti^2_{T,\de,2\eps}\cap\Ai_\ga\neq\emptyset\big]
\leq\sum_{k=k_-(\eps)}^{k_+(\eps)-1}\P\big[\Si^2_{T,\de,\eps,k\eps}\cap\Ai_\ga\neq\emptyset\big],
\ee
from which we see that (\ref{ctight}) implies (\ref{cS}).
\epro

\bpro[of Theorem~\ref{T:tight}]
By Theorem~\ref{T:Rtight} and Lemma~\ref{L:alt} it suffices to prove that (\ref{tight})~(i) is equivalent to
\be\label{MS2}
\lim_{\de\to 0}\sup_{\ga\in\Ga}\P\big[\Ai_\ga\cap\Ti^{\rm M}_{T,\de,\eta}\neq\emptyset\big]=0\quad\forall T,\eta>0
\ee
and (\ref{tight})~(ii) is equivalent to
\be\label{JS2}
\lim_{\de\to 0}\sup_{\ga\in\Ga}\P\big[\Ai_\ga\cap\Ti^{\rm J}_{T,\de,\eta}\neq\emptyset\big]=0\quad\forall T,\eta>0.
\ee
By formulas (\ref{setcom1}) and (\ref{setcom3}) of Lemma~\ref{L:setcom}, for each $T,\de,\eps>0$ and $r\in\R$, there exist $\eta,\eta'>0$ such that
\be
\P\big[\Si^{\rm M}_{T,\de,\eps,r}\cap\Ai_\ga\neq\emptyset\big]
\leq\P\big[\Ti^{\rm M}_{T,\de,\eta}\cap\Ai_\ga\neq\emptyset\big]
\quand
\P\big[\Si^{\rm J}_{T,\de,\eps,r}\cap\Ai_\ga\neq\emptyset\big]
\leq\P\big[\Ti^{\rm J}_{T,\de,\eta'}\cap\Ai_\ga\neq\emptyset\big],
\ee
from which we see that (\ref{MS2}) implies (\ref{tight})~(i) and (\ref{JS2}) implies (\ref{tight})~(ii). The converse implications follow from formulas (\ref{setcom2}) and (\ref{setcom4}) of Lemma~\ref{L:setcom} which tell us that
\be\ba{l}
\dis\P\big[\Ti^{\rm M}_{T,\de,\eta}\cap\Ai_\ga\neq\emptyset\big]\leq\sum_{k=k_-(\eps)}^{k_+(\eps)-1}\P\big[\Si^{\rm M}_{T,\de,\eps,k\eps}\cap\Ai_\ga\neq\emptyset\big],\\[5pt]
\dis\P\big[\Ti^{\rm J}_{T,\de,\eta}\cap\Ai_\ga\neq\emptyset\big]\leq\sum_{k=k_-(\eps)}^{k_+(\eps)-3}\P\big[\Si^{\rm J}_{T,\de,\eps,k\eps}\cap\Ai_\ga\neq\emptyset\big].
\ec
\epro

\section{Proof of the non-crossing criterion}\label{S:crosproof}

\subsection{The Fell topology}

Let $\Xc$ be a metrisable space that is separable and locally compact and let $\Xc^\ast:=\Xc\cup\{\infty\}$ denote its one-point compactification. It is well-known that $\Xc^\ast$ is metrisable. We observe that $\Ci:=\{A\in\Ki_+(\Xc^\ast):\infty\in A\}$ is a closed subset of $\Ki_+(\Xc^\ast)$ and hence compact. The map $A\mapsto A^\ast:=A\cup\{\infty\}$ is a bijection from the space $\Cl(\Xc)$ of closed subsets of $\Xc$ to $\Ci$. We use this to equip $\Cl(\Xc)$ with a metric $d_{\rm F}$ defined as
\be
d_{\rm F}(A,B):=d_{\rm H}\big(A^\ast,B^\ast\big)\qquad\big(A,B\in\Cl(\Xc)\big),
\ee
where $d_{\rm H}$ is the Hausdorff metric on $\Ki_+(\Xc^\ast)$ associated with an arbitrary metric $d$ generating the topology on $\Xc^\ast$. By our previous observation, $\big(\Cl(\Xc),d_{\rm F}\big)$ is compact. The \emph{Fell topology} is the weakest topology on $\Cl(\Xc)$ such that (i) $\{A:A\cap C=\emptyset\}$ is open for each compact $C\sub\Xc$ and (ii) $\{A:A\cap U\neq\emptyset\}$ is open for each open $U\sub\Xc$.

\bl[Fell topology]
The\label{L:Fell} metric $d_{\rm F}$ generates the Fell topology on $\Cl(\Xc)$.
\el

\bpro
If we equip $\Ki_+(\Xc^\ast)$ with the Hausdorff topology, then the topology generated by $d_{\rm F}$ is the initial topology for the map $\Cl(\Xc)\ni A\mapsto A^\ast\in\Ki_+(\Xc^\ast)$. Since $\Xc^\ast$ is compact, it is well-known that the Hausdorff and Fell topologies on $\Ki_+(\Xc^\ast)$ coincide, so the topology generated by $d_{\rm F}$ is the weakest topology on $\Cl(\Xc)$ such that (i) $\{A\in\Cl(\Xc):A^\ast\cap C=\emptyset\}$ is open for each compact $C\sub\Xc^\ast$ and (ii) $\{A\in\Cl(\Xc):A^\ast\cap U\neq\emptyset\}$ is open for each open $U\sub\Xc^\ast$. Conditions (i) and (ii) are trivially fulfilled if $C$ or $U$ contain the point $\infty$. On the other hand, the set of compact $C\sub\Xc^\ast$ for which $\infty\not\in C$ is precisely the set of compact subsets of $\Xc$, and the set of open $U\sub\Xc^\ast$ for which $\infty\not\in U$ is precisely the set of open subsets of $\Xc$, so the definition of the initial topology for the map $A\mapsto A^\ast$ gives us precisely the definition of the Fell topology on $\Cl(\Xc)$.
\epro

We have the following useful convergence criterion.

\bl[Convergence in the Fell topology]
Let\label{L:locHau} $A_n,A\in\Cl(\Xc)$. Then $A_n\to A$ in the Fell topology if and only if
\bc\label{locHau}
A&\sub&\dis\big\{x\in\Xc:\exists x_n\in A^\ast_n\mbox{ s.t.\ }x_n\to x\big\}\\[5pt]
&\sub&\dis\big\{x\in\Xc:\exists x_n\in A^\ast_n
\mbox{ s.t.\ $x$ is a cluster point of } (x_n)_{n\in\N}\big\}\sub A.
\ec
If $A_n\neq\emptyset$ for all $n$, then the statement remains true if in (\ref{locHau}) one replaces $A^\ast_n$ by $A_n$.
\el

\bpro
If (\ref{locHau}) holds, then
\bc\label{Ainf}
A^\ast&\sub&\dis\big\{x\in\Xc^\ast:\exists x_n\in A^\ast_n\mbox{ s.t.\ }x_n\to x\big\}\\[5pt]
&\sub&\dis\big\{x\in\Xc^\ast:\exists x_n\in A^\ast_n
\mbox{ s.t.\ $x$ is a cluster point of } (x_n)_{n\in\N}\big\}\sub A^\ast
\ec
and conversely (\ref{Ainf}) implies (\ref{locHau}), so the claim follows from Lemma~\ref{L:Hauconv} and the definition of $d_{\rm F}$. If $x_n\in\Xc^\ast$ satisfy $x_n\to x$ for some $x\in\Xc$, then $x_n\in\Xc$ for all $n$ large enough. Therefore, provided that $A_n\neq\emptyset$ for all $n$, in (\ref{locHau}) we may replace $A^\ast_n$ by $A_n$.
\epro

\subsection{Non-crossing paths}

In this subsection we prepare for the proof of Theorem~\ref{T:noncros}. For $\pi_1,\pi_2\in\Pi^|$, we say that $\pi_1$ \emph{crosses} $\pi_2$ if $\pi_1\not\lef\pi_2$ and $\pi_2\not\lef\pi_1$, that is, if the set $\{\pi_1,\pi_2\}$ is not non-crossing in the sense defined in Subsection~\ref{S:noncros}. We say that $\pi_1$ \emph{collides} with $\pi_2$ \emph{at time} $t$ if $t\in I_{\pi_1}\cap I_{\pi_2}$ and $\pi_1(t\pm),\pi_2(t\mp)<\pi_1(t\mp),\pi_2(t\pm)$, where the sign $\pm$ can be either $+$ or $-$ and $\mp$ is the opposite sign. In words, this says that at time $t$ the paths $\pi_1$ and $\pi_2$ jump over some interval in opposite directions. It is not hard to see that if $\pi_1$ collides with $\pi_2$ at some time $t$, then $\pi_1$ crosses $\pi_2$. We let $I^\circ$ denote the interior of a closed real interval and let $\pa I:=I\beh I^\circ$ denote the set of its finite boundary points. We will prove the following result. Note that this shows in particular that for bi-infinite paths, the non-crossing property is preserved under limits.

\bp[Crossing in the limit]
Assume\label{P:collis} that $\pi^n_1,\pi^n_2\in\Pi^|$ satisfy $\pi^n_i\to\pi_i$ as $n\to\infty$ $(i=1,2)$ in the M1 topology for some $\pi_1,\pi_2\in\Pi^|$, and that $\pi^n_1$ does not cross $\pi^n_2$ for any $n$. Then precisely one of the following statements must hold:
\begin{enumerate}
\item $\pi_1$ does not cross $\pi_2$,
\item $\pi_1$ collides with $\pi_2$ at some time $t\in\pa I_{\pi_1}\cup\pa I_{\pi_2}$.
\end{enumerate}
\ep

Proposition~\ref{P:collis} extends \cite[Lemma~3.4.9]{FS24}. It is not hard to guess the statement of Proposition~\ref{P:collis} after scribbling a few figures on a piece of paper, but giving a precise proof is a bit more work. Our proof follows a clear strategy but is nevertheless quite long, so we challenge the reader to find a shorter one.

\begin{figure}[t]
\begin{center}
\includegraphics{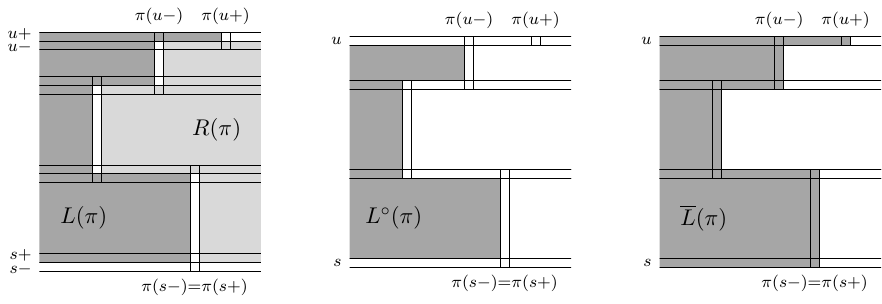}
\caption{Schematic depiction of the sets $L(\pi)$ and $R(\pi)$ from (\ref{LRdef}) and the sets $L^\circ(\pi)$ and $\ov L(\pi)$ from (\ref{LL}) for a piecewise constant path $\pi$. In the pictures for $L^\circ(\pi)$ and $\ov L(\pi)$, values attained by $\pi$ and times when $\pi$ makes a jump are enlarged to indicate whether they are included in the sets $L^\circ(\pi)$ and $\ov L(\pi)$. In the picture for $L(\pi)$ and $R(\pi)$, times $t$ when $\pi$ makes a jump are represented by two points, $t-$ and $t+$.}
\label{fig:LR}
\end{center}
\end{figure}

To prepare for the proof of Proposition~\ref{P:collis}, we start by giving a more direct characterisation of the relation $\pi_1\lef\pi_2$. The \emph{split real line} is the set $\R_\mfs$ consisting of all words of the form $t\star$ with $t\in\R$ and $\star\in\{-,+\}$. For each $\pi\in\Pi^|$ with starting time $s:=\inf\hat I_\pi$ and final time $u:=\sup\hat I_\pi$, we define $I^\mfs_\pi\sub\R_\mfs$ by
\be\label{Is}
I^\mfs_\pi:=\big\{t\star:t\in I_\pi,\ \star\in\{-,+\},\ t\star\neq s-,u+\big\},
\ee
and we set
\be
\label{Ilr}
\begin{alignedat}{4}
\dis I^{\rm l}_\pi
&
:=\dis\big\{t\star\,:\;
&&
t\star\in I^\mfs_\pi,
&&
\quad\text{or}\quad t\star=s-\text{ and }\pi(s+)<\pi(s-),\\
& && &&\quad\text{or}\quad t\star=u+\text{ and }\pi(u-)<\pi(u+)\big\},\\[8pt]
\dis I^{\rm r}_\pi
&
:=\dis\big\{t\star\,:\,
&&
\dis t\star\in I^\mfs_\pi,
&&
\quad\text{or}\quad t\star=s-\text{ and }\pi(s-)<\pi(s+),\\
& && 
&&\dis\quad\text{or}\quad t\star=u+\text{ and }\pi(u+)<\pi(u-)\big\},
\end{alignedat}
\ee
that is, we include $s-$ in $I^{\rm l}_\pi$ only if $s\in\R$ and $\pi(s+)<\pi(s-)$, and so on. Finally, we set
\bc\label{LRdef}
\dis L(\pi)&:=&\dis\big\{(x,t\pm)\in\ov\R\times I^{\rm l}_\pi:x<\pi(t\pm)\big\},\\[5pt]
\dis R(\pi)&:=&\dis\big\{(x,t\pm)\in\ov\R\times I^{\rm r}_\pi:\pi(t\pm)<x\big\}.
\ec
These definitions are illustrated in Figure~\ref{fig:LR}. Note that they simplify a lot for bi-infinite paths. The following lemma gives a more direct characterisation of the relation $\pi_1\lef\pi_2$.

\bl[Ordering of paths]
Two\label{L:piord} paths $\pi_1,\pi_2\in\Pi^|$ satisfy $\pi_1\lef\pi_2$ if and only if $L(\pi_1)\cap R(\pi_2)=\emptyset$.
\el

\bpro
The proof is only a slight extension of \cite[Lemma~3.3.2]{FS24} but for completeness we include a formal proof. For bi-infinite paths $\pi_1,\pi_2\in\Pi^\updo$ the statement is trivial. It is also straightforward to check that if $\pi'$ extends $\pi$, then $L(\pi)\sub L(\pi')$ and $R(\pi)\sub R(\pi')$. In view of this, the condition $L(\pi_1)\cap R(\pi_2)=\emptyset$ is clearly necessary for $\pi_1\lef\pi_2$. To prove sufficiency, we will prove the following statement:
\begin{itemize}
\item[(!)] Assume that $\pi_1,\pi_2\in\Pi^|$ satisfy $L(\pi_1)\cap R(\pi_2)=\emptyset$. Then there exists an extension $\pi'_1$ of $\pi_1$ such that $\pi'_1\in\Pi^\up$ and $L(\pi'_1)\cap R(\pi_2)=\emptyset$. 
\end{itemize}
By symmetry, it then follows that we can also extend $\pi_2$ so that it runs till time $+\infty$ and in two next steps we can extend $\pi_1$ and $\pi_2$ so that they both start at time $-\infty$, all the time while preserving the property that $L(\pi_1)\cap R(\pi_2)=\emptyset$, which for bi-infinite paths trivially implies $\pi_1\lef\pi_2$. It therefore remains to prove (!).

Let $s_i:=\inf I_{\pi_i}$ and $u_i:=\sup I_{\pi_i}$. The cases $u_1=\pm\infty$ are trivial so without loss of generality we assume that $u_1\in\R$. The cases $u_1<s_2$ and $u_1>u_2$ are also trivial so without loss of generality we can assume that $s_2\leq u_1\leq u_2$. If $\pi_1(u_1+)\leq\pi_1(u_1-)$, then we can simply extend $\pi_1$ to the path that jumps to $-\infty$ at time $u_1$ and stays there, so we are done. We therefore from now on assume $\pi_1(u_1-)<\pi_1(u_1+)$. We distinguish two cases: I.\ $u_1=u_2$, and II.\ $u_1<u_2$. In case~I we can extend $\pi_1$ to the path that is constantly equal to $\pi_1(u_1+)$ after time $u_1$. Since $\pi_1(u_1-)<\pi_1(u_1+)$, extending $\pi_1$ in this way will only add points to the set $L(\pi_1)$ whose time coordinates are strictly larger than $u_1$ and that hence will not lie inside $R(\pi_2)$. In case~II the condition $L(\pi_1)\cap R(\pi_2)=\emptyset$ implies $\pi_1(u_1+)\leq\pi_2(u_1+)$, so we can extend $\pi_1$ to the path that jumps to $\pi_2(u_1+)$ at time $u_1$ and then stays equal to $\pi_2$ until time $u_2+$, reducing case~II to case~I.
\epro

%
%

We next start to investigate how crossing behaves under limits. For any path $\pi$, we let $I^\circ_\pi$ denote the interior of its domain $I_\pi\sub\R$. To simplify notation, we set $\Sp:=\ov\R\times\R$. For any $\pi\in\Pi^|$, we write
\bc\label{LL}
\dis L^\circ(\pi)&:=&\dis\big\{(x,t)\in\Sp:t\in I^\circ_\pi,\ x<\pi(t-)\wedge\pi(t+)\big\},\quad L^{\rm c}(\pi):=\Sp\beh L^\circ(\pi),\\[5pt]
\dis\ov L(\pi)&:=&\dis\big\{(x,t)\in\Sp:t\in I_\pi,\ x\leq\pi(t-)\vee\pi(t+)\big\}.
\ec
These definitions are illustrated in Figure~\ref{fig:LR}. We define $R^\circ(\pi),R^{\rm c}(\pi)$, and $\ov R(\pi)$ analogously, with $x<\pi(t-)\wedge\pi(t+)$ replaced by $\pi(t-)\vee\pi(t+)<x$ and so forth. They look a lot like the sets $L(\pi)$ and $R(\pi)$ from (\ref{LRdef}), but as we will see in Lemma~\ref{L:leftconv} below, they are better behaved under limits. Note that $L^\circ(\pi)$ is open while $L^{\rm c}(\pi)$ and $\ov L(\pi)$ are closed.

\bl[Consequence of ordering]
If\label{L:intcros} $\pi_1,\pi_2\in\Pi^|$ satisfy $\pi_1\lef\pi_2$, then
\be
L^\circ(\pi_1)\cap\ov R(\pi_2)=\emptyset
\quand
\ov L(\pi_1)\cap R^\circ(\pi_2)=\emptyset.
\ee
\el

\bpro
For any path $\pi\in\Pi^|$ we set $L^\pm(\pi):=\{(x,t):(x,t\pm)\in L(\pi)\}$ and we define $R^-(\pi)$ and $R^+(\pi)$ similarly. It is straightforward to check (see Figure~\ref{fig:LR}) that
\be\label{intclos}
L^\circ(\pi)\mbox{ is the interior of }L^-(\pi)\cap L^+(\pi)
\quand
\ov L(\pi)\mbox{ is the closure of }L^-(\pi)\cup L^+(\pi),
\ee
and likewise for $R^\circ(\pi)$ and $\ov R(\pi)$. The condition $L(\pi_1)\cap R(\pi_2)=\emptyset$ is equivalent to $L^-(\pi_1)\cap R^-(\pi_2)=\emptyset$ and $L^+(\pi_1)\cap R^+(\pi_2)=\emptyset$, which implies
\be
\big(L^-(\pi_1)\cap L^+(\pi_1)\big)\cap\big(R^-(\pi_2)\cup R^+(\pi_2)\big)=\emptyset.
\ee
Combining this with (\ref{intclos}) we see that $\pi_1\lef\pi_2$ implies $L^\circ(\pi_1)\cap\ov R(\pi_2)=\emptyset$. The proof that $\ov L(\pi_1)\cap R^\circ(\pi_2)=\emptyset$ is the same.
\epro


\bl[Areas left of a path]
Assume\label{L:leftconv} that $\pi_n\in\Pi^|$ satisfy $\pi_n\to\pi$ as $n\to\infty$ in the M1 topology for some $\pi\in\Pi^|$. Then $\ov L(\pi_n)\to\ov L(\pi)$ and $L^{\rm c}(\pi_n)\to L^{\rm c}(\pi)$ in the Fell topology on $\Cl(\Sp)$.
\el

\bpro
Assume that $\pi_n\to\pi$ in the M1 topology. Then $\ov\pi_n\to\ov\pi$ in the Hausdorff topology by Lemma~\ref{L:Haupi}, which by Lemma~\ref{L:Hauconv} means that:
\begin{itemize}
\item[(i)] for each $z\in\ov\pi$ there exist $z_n\in\ov\pi_n$ such that $z_n\to z$,
\item[(ii)] if $z_n\in\ov\pi_n$ have a cluster point $z$, then $z\in\ov\pi$.
\end{itemize}
By Lemma~\ref{L:contim} moreover $\hat I_{\pi_n}\to\hat I_\pi$ in the Hausdorff topology on $\Ki_+(\ov\R)$.

We will only prove that $L^{\rm c}(\pi_n)\to L^{\rm c}(\pi)$. The proof that $\ov L(\pi_n)\to\ov L(\pi)$ is similar, but easier. We note that
\be
L^{\rm c}(\pi)=\big\{(x,t)\in\Sp:t\not\in I^\circ_\pi\big\}\cup\big\{(x,t)\in\Sp:t\in I^\circ_\pi\mbox{ and }\pi(t-)\wedge\pi(t+)\leq x\big\}.
\ee
By Lemma~\ref{L:locHau}, we have to show that:
\begin{itemize}
\item[(i)'] for each $z\in L^{\rm c}(\pi)$ there exist $z_n\in L^{\rm c}(\pi_n)^\ast$ such that $z_n\to z$,
\item[(ii)'] if $z_n\in L^{\rm c}(\pi_n)^\ast$ have a cluster point $z\in\Sp$, then $z\in L^{\rm c}(\pi)$.
\end{itemize}
To prove (i)', assume that $(x,t)\in L^{\rm c}(\pi)$. Then either I.\ $t\not\in I^\circ_\pi$ or II.\ $t\in I^\circ_\pi$ and $\pi(t-)\wedge\pi(t+)\leq x$. In case~I, we can use that $\hat I_{\pi_n}\to\hat I_\pi$ to choose $t_n\not\in I^\circ_{\pi_n}$ such that $t_n\to t$. Now $(x,t_n)\in L^{\rm c}(\pi_n)$ satisfy $(x,t_n)\to(x,t)$. In case~II, we set $y:=\pi(t-)\wedge\pi(t+)$. Then $(y,t)\in\ov\pi$ and hence by (i) there exist $(y_n,t_n)\in\ov\pi_n$ such that $(y_n,t_n)\to(y,t)$. Since $y\leq x$ it follows that $L^{\rm c}(\pi_n)\ni (y_n\vee x,t_n)\to(x,t)$.

To prove (ii)', assume that $(x_n,t_n)\in L^{\rm c}(\pi_n)$ have a cluster point $(x,t)$. By going to a subsequence, we may assume that $(x_n,t_n)\to(x,t)$ and either I.\ $t_n\not\in I^\circ_{\pi_n}$ for all $n$ or II.\ $t_n\in I^\circ_{\pi_n}$ and $\pi_n(t_n-)\wedge\pi_n(t_n+)\leq x_n$ for all $n$. In case~I by the fact that $\hat I_{\pi_n}\to\hat I_\pi$ we have $t=\lim_{n\to\infty}t_n\not\in I^\circ_\pi$ so $(x,t)\in L^{\rm c}(\pi)$ and we are done. Case~II is trivial if $t\not\in I^\circ_\pi$ so without loss of generality we assume $t\in I^\circ_\pi$. Then, by going to a further subsequence if necessary, we can assume that $y_n:=\pi_n(t_n-)\wedge\pi_n(t_n+)$ converge to a limit $y\in\ov\R$. Then $\ov\pi_n\ni(y_n,t_n)\to(y,t)$ so by (i) $(y,t)\in\ov\pi$ and hence $\pi(t-)\wedge\pi(t+)\leq y$. Since $y_n\leq x_n$ we have $y\leq x$ so we conclude that $(x,t)\in L^{\rm c}(\pi)$ and we are done.
\epro

As a simple corollary of Lemma~\ref{L:leftconv}, we obtain:

\bcor[Limits of ordered paths]
Assume\label{C:limlef} that $\pi^n_1,\pi^n_2\in\Pi^|$ satisfy $\pi^n_1\lef\pi^n_2$ for all $n$ and that $\pi^n_i\to\pi_i$ as $n\to\infty$ $(i=1,2)$ in the M1 topology for some $\pi_1,\pi_2\in\Pi^|$. Then $L^\circ(\pi_1)\cap\ov R(\pi_2)=\emptyset$ and $\ov L(\pi_1)\cap R^\circ(\pi_2)=\emptyset$.
\ecor

\bpro
By symmetry, it suffices to prove that $L^\circ(\pi_1)\cap\ov R(\pi_2)=\emptyset$ or equivalently $\ov R(\pi_2)\sub L^{\rm c}(\pi_1)$. Since $\pi^n_1\lef\pi^n_2$ for all $n$, Lemma~\ref{L:intcros} tells us that $\ov R(\pi^n_2)\sub L^{\rm c}(\pi^n_1)$ for all $n$, so the claim follows from Lemma~\ref{L:leftconv}, using Lemma~\ref{L:locHau}.
\epro

\bpro[of Proposition~\ref{P:collis}]
Assume that $\pi^n_1,\pi^n_2\in\Pi^|$ satisfy $\pi^n_i\to\pi_i$ as $n\to\infty$ $(i=1,2)$ in the M1 topology for some $\pi_1,\pi_2\in\Pi^|$, and that $\pi^n_1$ does not cross $\pi^n_2$ for any $n$. It is clear from Lemma~\ref{L:piord} that (i) and (ii) cannot hold simultaneously, so it suffices to prove that if $\pi_1$ crosses $\pi_2$, then $\pi_1$ collides with $\pi_2$ at some time $t\in\pa I_{\pi_1}\cup\pa I_{\pi_1}$.

\begin{figure}[t]
\begin{center}
\includegraphics{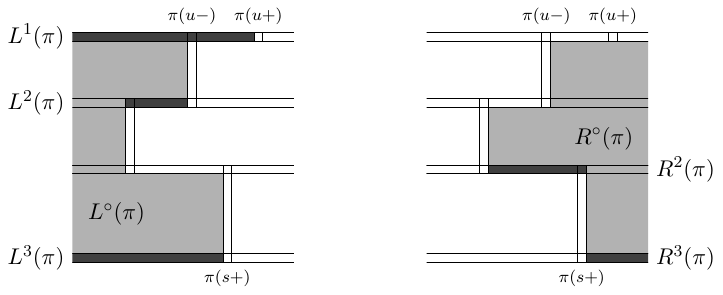}
\caption{Schematic depiction of the sets $L^1(\pi)$, $L^2(\pi)$, and $L^3(\pi)$ defined in (\ref{LLL}), and the analogously defined sets $R^1(\pi)$, $R^2(\pi)$, and $R^3(\pi)$. In this example $R^1(\pi)=\emptyset$ because $\pi(u-)\leq\pi(u+)$.}
\label{fig:LRplus}
\end{center}
\end{figure}

Since $\pi^n_1$ does not cross $\pi^n_2$ for any $n$, by going to a subsequence, we can assume that either $\pi^n_1\lef\pi^n_2$ for all $n$ or $\pi^n_2\lef\pi^n_1$ for all $n$. By symmetry between left and right, we can without loss of generality assume that we are in the first case. Then Corollary~\ref{C:limlef} tells us that
\be\label{cirov}
L^\circ(\pi_1)\cap\ov R(\pi_2)=\emptyset
\quand
\ov L(\pi_1)\cap R^\circ(\pi_2)=\emptyset.
\ee
Since $\pi_1$ crosses $\pi_2$, we have $\pi_1\not\lef\pi_2$ and $\pi_2\not\lef\pi_1$. For any path $\pi\in\Pi^|$ we define $L^\pm(\pi)$ and $R^\pm(\pi)$ as in the proof of Lemma~\ref{L:intcros}. Since $\pi_1\not\lef\pi_2$ we have $L^-(\pi_1)\cap R^-(\pi_2)\neq\emptyset$ or $L^+(\pi_1)\cap R^+(\pi_2)\neq\emptyset$. By symmetry with respect to time reversal, we may assume that we are in the second case. It is easy to see that $L^\circ(\pi_1)\sub L^+(\pi_1)\sub\ov L(\pi_1)$ and $R^\circ(\pi_1)\sub R^+(\pi_1)\sub\ov R(\pi_1)$ (see Figure~\ref{fig:LR}). Therefore the fact that $L^+(\pi_1)\cap R^+(\pi_2)\neq\emptyset$ and (\ref{cirov}) imply that
\be\label{difdif}
\big(L^+(\pi_1)\beh L^\circ(\pi_1)\big)\cap\big(R^+(\pi_1)\beh R^\circ(\pi_1)\big)\neq\emptyset.
\ee
As illustrated in Figure~\ref{fig:LRplus}, for any path $\pi$ with starting time $s:=\inf I_\pi$ and final time $u:=\sup I_\pi$, we can write $L^+(\pi)\beh L^\circ(\pi)=L^1(\pi)\cup L^2(\pi)\cup L^3(\pi)$, where
\bc\label{LLL}
\dis L^1(\pi)&:=&\dis\big\{(x,u):x<\pi(u+),\ \pi(u-)<\pi(u+)\big\},\\[5pt]
\dis L^2(\pi)&:=&\dis\big\{(x,t):s<t<u,\ \pi(t-)<x<\pi(t+),\ \pi(t-)<\pi(t+)\big\},\\[5pt]
\dis L^3(\pi)&:=&\dis\big\{(x,s):x<\pi(s+)\big\}.
\ec
Similarly, we can decompose $R^+(\pi)\beh R^\circ(\pi)$ into three sets $R^1(\pi)$, $R^2(\pi)$, and $R^3(\pi)$ which contain points $(x,t)$ with $t=u$, $s<t<u$, and $t=s$, respectively. 

We observe that for each $(x,t)\in L^2(\pi)\cup L^3(\pi)$, one has $(x,t+\eps)\in L^\circ(\pi)$ for all $\eps>0$ small enough. A similar statement holds for $(x,t)\in R^2(\pi)\cup R^3(\pi)$. Since $L^\circ(\pi_1)\cap R^\circ(\pi_1)=\emptyset$ by (\ref{cirov}), it follows that $(L^2(\pi_1)\cup L^3(\pi_1))\cap(R^2(\pi_1)\cup R^3(\pi_1))=\emptyset$. This allows us to improve (\ref{difdif}) to the more precise statement
\be\label{123}
L^1(\pi_1)\cap\big(R^1(\pi_2)\cup R^2(\pi_2)\cup R^3(\pi_2)\big)\neq\emptyset
\quad\mbox{or}\quad
\big(L^1(\pi_1)\cup L^2(\pi_1)\cup L^3(\pi_1)\big)\cap R^1(\pi_2)\neq\emptyset.
\ee
If $L^1(\pi_1)\cap R^3(\pi_2)\neq\emptyset$, then $\pi_2\lef\pi_1$, which contradicts the fact that $\pi_2\not\lef\pi_1$. The same argument gives $L^3(\pi_1)\cap R^1(\pi_2)=\emptyset$, so we can improve (\ref{123}) to
\be\label{12}
L^1(\pi_1)\cap\big(R^1(\pi_2)\cup R^2(\pi_2)\big)\neq\emptyset
\quad\mbox{or}\quad
\big(L^1(\pi_1)\cup L^2(\pi_1)\big)\cap R^1(\pi_2)\neq\emptyset,
\ee
which implies that $\pi_1$ collides with $\pi_2$ at some time $t\in\pa I_{\pi_1}\cup\pa I_{\pi_1}$.
\epro


\subsection{Tightness of non-crossing sets}

In this subsection we prove Theorem~\ref{T:noncros} and Lemma~\ref{L:biclos}. We start by showing that (\ref{crotig}) implies (\ref{tight})~(i).

\bl[Jumps in opposite directions]
Let\label{L:CS} $\Ai$ be a compact subset of $\Pi^|$, let $T,\de,\eps>0$ and $r\in\R$. Then
\be
\Si^{\rm M}_{T,\de,\eps,r}\cap\Ai\neq\emptyset
\quad\mbox{implies}\quad
\Ci^{\rm M}_{T,\de,\eps,r}\cap(\Ai\times\Ai)\neq\emptyset.
\ee
\el

\bpro
One has $\Si^{+-}_{T,\de,\eps,r}\cap\Ai\neq\emptyset$ if and only if:
\begin{itemize}
\item[(S)] There exists a $\pi\in\Ai$ and $(x,s),(y,t),(z,u)\in\pi$ with $(x,s)\pre(y,t)\pre(z,u)$, $-T\leq s\leq t\leq u\leq T$, and $u-s\leq\de$, such that $x,z\leq r$ and $r+\eps\leq y$. 
\end{itemize}
Similarly, one has $\Ci^{\rm M}_{T,\de,\eps,r}\cap(\Ai\times\Ai)\neq\emptyset$ if and only if:
\begin{itemize}
\item[(C)] For $i=1,2$ there exists $\pi_i\in\Ai$ and $(x_i,s_i),(y_i,t_i)\in\pi_i$ with $(x_i,s_i)\pre(y_i,t_i)$, $-T\leq s_i\leq t_i\leq T$, and $(t_1\vee t_2)-(s_1\wedge s_2)\leq\de$, such that $x_1,y_2\leq r$ and $r+\eps\leq y_1,x_2$. 
\end{itemize}
If (S) holds, then setting $\pi_1=\pi_2:=\pi$, $(x_1,s_1):=(x,s)$, $(y_1,t_1)=(x_2,s_2):=(y,t)$, and $(y_2,t_2):=(z,u)$, we see that (C) holds. The same argument with the roles of $\pi_1$ and $\pi_2$ interchanged shows that $\Si^{-+}_{T,\de,\eps,r}\cap\Ai\neq\emptyset$ implies (C).
\epro

\bl[Crossing in the limit]
Let\label{L:limcros} $\Pi^|$ be equipped with the M1 topology and let $\Ai_n$ be non-crossing compact subsets of $\Pi^|$. Assume that $\Ai_n\to\Ai$ in the Hausdorff topology on $\Ki_+(\Pi)$. Then $\Ai\sub\Pi^|$. If $\Ai$ is non-crossing, then
\be\label{limnoncros}
\forall T,\eps>0,\ r\in\R\ \exists\de>0,\ m\in\N\mbox{ s.t.\ }\Ci^{\rm M}_{T,\de,\eps,r}\cap(\Ai_n\times\Ai_n)=\emptyset\ \forall n\geq m.
\ee
On the other hand, if $\Ai$ is not non-crossing, then
\be\label{limcros}
\exists T,\eps>0,\ r\in\R\mbox{ s.t.\ } \forall\de>0\ \exists m\in\N\mbox{ s.t.\ }\Ci^{\rm M}_{T,\de,\eps,r}\cap(\Ai_n\times\Ai_n)\neq\emptyset\ \forall n\geq m.
\ee
\el

\bpro
It is easy to see that $\Pi^|$ is a closed subset of $\Pi$. As a result, by grace of Lemma~\ref{L:Hauconv}, $\Ki_+(\Pi^|)$ is a closed subset of $\Ki_+(\Pi)$. This proves that $\Ai\sub\Pi^|$.

Let us say that $\Ai$ \emph{contains colliding paths} if there exist $\pi_1,\pi_2\in\Ai$ and $t\in\R$ such that $\pi_1$ collides with $\pi_2$ at time $t$. We will prove the following statements:
\begin{itemize}
\item[I.] $\Ai$ is non-crossing if and only if $\Ai$ does not contain colliding paths.
\item[II.] If $\Ai$ does not contain colliding paths, then (\ref{limnoncros}) holds.
\item[III.] If $\Ai$ contains colliding paths, then (\ref{limcros}) holds.
\end{itemize}
We start by proving I. It is clear that if $\Ai$ contains colliding paths, then $\Ai$ is not non-crossing. To prove the converse, assume that $\Ai$ does not contain colliding paths and let $\pi_1,\pi_2\in\Ai$. By Lemma~\ref{L:Hauconv}, there exist $\pi_1^n,\pi_2^n\in\Ai_n$ such that $\pi_i^n\to\pi_i$ $(i=1,2)$. Since $\Ai_n$ is non-crossing, for each $n$, the path $\pi_1^n$ does not cross $\pi_2^n$. By Proposition~\ref{P:collis}, it follows that either $\pi_1$ collides with $\pi_2$ or $\pi_1$ does not cross $\pi_2$. The first option is excluded by our assumption that $\Ai$ does not contain colliding paths, so we conclude that $\Ai$ is non-crossing.

We next prove II. We will show that if (\ref{limnoncros}) does not hold, then $\Ai$ contains colliding paths. If (\ref{limnoncros}) does not hold, then there exist $T,\eps>0$ and $r\in\R$ such that
\be
\forall\de>0,\ m\in\N\ \exists n\geq m\mbox{ s.t.\ }\Ci^{\rm M}_{T,\de,\eps,r}\cap(\Ai_n\times\Ai_n)\neq\emptyset.
\ee
Fix $\de_k>0$ such that $\de_k\to 0$. Then we can choose $n(k)\to\infty$ such that $\Ci^{\rm M}_{T,\de_n,\eps,r}\cap(\Ai_{n(k)}\times\Ai_{n(k)})\neq\emptyset$ for each $k$. This implies that for $i=1,2$ there exists $\pi^k_i\in\Ai_{n(k)}$ and $(x^k_i,s^k_i),(y^k_i,t^k_i)\in\pi^k_i$ with $(x^k_i,s^k_i)\pre(y^k_i,t^k_i)$, $-T\leq s^k_i\leq t^k_i\leq T$, and $(t^k_1\vee t^k_2)-(s^k_1\wedge s^k_2)\leq\de_k$, such that $x^k_1,y^k_2\leq r$ and $r+\eps\leq y^k_1,x^k_2$. By Lemma~\ref{L:Hauconv} and the fact that $\Ai_{n(k)}\to\Ai$, there exists a compact $\Ci\sub\Pi$ such that $\Ai_{n(k)}\sub\Ci$ for all $k$. In view of this, using also the compactness of $\ov\R$ and $[-T,T]$, by going to a subsequence if necessary, we can assume that as $k\to\infty$
\be
\pi^k_i\to\pi_i,\ (x^k_i,s^k_i)\to(x_i,s_i),\quad(y^k_i,t^k_i)\to(y_i,t_i)\quad(i=1,2).
\ee
By Lemma~\ref{L:Hauconv} we have $\pi_1,\pi_2\in\Ai$. Combining Lemma~\ref{L:Hauconv} with Lemma~\ref{L:JMconv} we see that $(x_i,s_i),(y_i,t_i)\in\ov\pi_i$ with $(x_i,s_i)\pre(y_i,t_i)$. Taking the limit in $x^k_1,y^k_2\leq r$ and $r+\eps\leq y^k_1,x^k_2$ we obtain that $x_1,y_2\leq r$ and $r+\eps\leq y_1,x_2$. Finally, since $(t^k_1\vee t^k_2)-(s^k_1\wedge s^k_2)\leq\de_k$, we must have $s_1=s_2=t_1=t_2=:t$ for some $t\in[-T,T]$. We have to be a bit careful since as a result of using the M1 topology we only know that $(x_i,s_i),(y_i,t_i)$ are elements of the filled graph $\ov\pi_i$, and not necessarily of the closed graph $\pi_i$. Nevertless, the properties we have proved are enough to conclude that $\pi_1$ collides with $\pi_2$ at time $t$.

It remains to prove III. Assume that $\pi_1,\pi_2\in\Ai$ collide at time $t$. Then $t\in[-T/2,T/2]$ for $T$ large enough and (possibly after interchanging the roles of $\pi_1$ and $\pi_2$) there exist $x_i,y_i$ with $(x_i,t),(y_i,t)\in\ov\pi_i$ and $(x_i,t)\pre(y_i,t)$ $(i=1,2)$ such that $x_1,y_2\leq r-\eps$ and $r+2\eps\leq y_1,x_2$ for some $\eps>0$ and $r\in\R$. By Lemma~\ref{L:Hauconv} there exist $\pi^n_1,\pi^n_2\in\Ai_n$ such that $\pi^n_i\to\pi_i$ $(i=1,2)$, so by Lemmas \ref{L:Hauconv} and \ref{L:JMconv} there exist $(x^n_i,s^n_i),(y^n_i,t^n_i)\in\ov\pi^n_i$ with $(x^n_i,s^n_i)\pre(y^n_i,t^n_i)$ such that $x^n_i\to x_i$ and $y^n_i\to y_i$ $(i=1,2)$. Then for each $\de>0$ one has $(t^n_1\vee t^n_2)-(s^n_1\wedge s^n_2)\leq\de$ for $n$ large enough. Moreover, for all $n$ large enough $s^n_i,t^n_i\in[-T,T]$, $x^n_1,y^n_2\leq r$, and $r+\eps\leq y^n_1,x^n_2$. By making $x^n_1,y^n_2$ smaller if necessary and $y^n_1,x^n_2$ larger if necessary, we can make sure that $(x^n_i,s^n_i),(y^n_i,t^n_i)$ are elements of the closed graph $\pi^n_i$ and not just of the filled graph $\ov\pi^n_i$, while preserving all other properties mentioned above. Then
\be
\forall\de>0\ \exists m\in\N\mbox{ s.t.\ }\Ci^{\rm M}_{T,\de,\eps,r}\cap(\Ai_n\times\Ai_n)\neq\emptyset\ \forall n\geq m,
\ee
which shows that (\ref{limcros}) holds.
\epro

\bpro[of Theorem~\ref{T:noncros}]
Set $\Mi:=\{\mu_\ga:\ga\in\Ga\}$. Let $\Mi_1(\Ki_+(\Pi))$ denote the space of probability measures on $\Ki_+(\Pi)$, equipped with the topology of weak convergence. We naturally view $\Mi_1(\Ki_{\rm nc}(\Pi^|))$ as a subset of $\Mi_1(\Ki_+(\Pi))$. Then $\Mi$ is precompact as a subset of $\Mi_1(\Ki_{\rm nc}(\Pi^|))$ if and only if $\Mi$ is precompact as a subset of $\Mi_1(\Ki_+(\Pi))$ and the closure of $\Mi$ is contained in $\Mi_1(\Ki_{\rm nc}(\Pi^|))$. Thus, by Prokhorov's theorem, $\Mi$ is tight on $\Ki_{\rm nc}(\Pi^|)$ if and only if $\Mi$ is tight on $\Ki_+(\Pi)$ and the closure of $\Mi$ is concentrated on $\Ki_{\rm nc}(\Pi^|)$. Therefore, by Theorem~\ref{T:tight}, it suffices to prove the following statements:
\begin{itemize}
\item[I.] If (\ref{crotig}) holds, then (\ref{tight})~(i) holds and the closure of $\Mi$ is concentrated on $\Ki_{\rm nc}(\Pi^|)$.
\item[II.] If (\ref{tight})~(i) holds and the closure of $\Mi$ is concentrated on $\Ki_{\rm nc}(\Pi^|)$, then (\ref{crotig}) holds.
\end{itemize}
We first prove I. Assume that (\ref{crotig}) holds. Then by Lemma~\ref{L:CS} condition (\ref{tight})~(i) is satisfied, so it remains to show that the closure of $\Mi$ is concentrated on $\Ki_{\rm nc}(\Pi^|)$. Assume that $\mu_n\in\Mi$ converge weakly to a probability law $\mu$ on $\Ki_+(\Pi)$. By Skorokhod's representation theorem \cite[Cor~3.1.6 and Thm~3.1.8]{EK86}, we can couple random variables $\Ai_n,\Ai$ with laws $\mu_n,\mu$ such that $\Ai_n\to\Ai$ a.s. We need to show that $\Ai\in\Ki_{\rm nc}(\Pi^|)$ a.s. By Lemma~\ref{L:limcros} we have $\Ai\sub\Pi^|$ a.s. Therefore, by condition (\ref{limcros}) of Lemma~\ref{L:limcros}, to prove that $\Ai\in\Ki_{\rm nc}(\Pi^|)$ a.s., it suffices to show that
\be
\P\big[\exists T,\eps>0,\ r\in\R\mbox{ s.t.\ }\forall\de>0,\ \exists m\in\N\mbox{ s.t.\ }\Ci^{\rm M}_{T,\de,\eps,r}\cap(\Ai_n\times\Ai_n)\neq\emptyset\ \forall n\geq m\big]=0.
\ee
It suffices to check the condition (\ref{limcros}) for countably many values of $T,\eps$, and $r$ only: in particular, we can take $T=N$, $\eps=1/n$, and $r=k/(3n)$ with $N,n$ positive integers and $k\in\Z$. In view of this, it suffices to show that
\be
\P\big[\forall\de>0,\ \exists m\in\N\mbox{ s.t.\ }\Ci^{\rm M}_{T,\de,\eps,r}\cap(\Ai_n\times\Ai_n)\neq\emptyset\ \forall n\geq m\big]=0\quad\forall T,\eps>0,\ r\in\R.
\ee
Since $\de'\leq\de$ implies $\Ci^{\rm M}_{T,\de',\eps,r}\sub\Ci^{\rm M}_{T,\de,\eps,r}$, we can rewrite our previous formula as
\be\label{Plimcros}
\lim_{\de\to 0}\P\big[\exists m\in\N\mbox{ s.t.\ }\Ci^{\rm M}_{T,\de,\eps,r}\cap(\Ai_n\times\Ai_n)\neq\emptyset\ \forall n\geq m\big]=0\quad\forall T,\eps>0,\ r\in\R.
\ee
To see that (\ref{crotig}) implies (\ref{Plimcros}), we estimate
\be\ba{l}
\dis\P\big[\exists m\in\N\mbox{ s.t.\ }\Ci^{\rm M}_{T,\de,\eps,r}\cap(\Ai_n\times\Ai_n)\neq\emptyset\ \forall n\geq m\big]=\lim_{m\to\infty}\P\big[\Ci^{\rm M}_{T,\de,\eps,r}\cap(\Ai_n\times\Ai_n)\neq\emptyset\ \forall n\geq m\big]\\[5pt]
\dis\quad\leq\limsup_{m\to\infty}\P\big[\Ci^{\rm M}_{T,\de,\eps,r}\cap(\Ai_m\times\Ai_m)\neq\emptyset\big]\leq\sup_{m\in\N}\P\big[\Ci^{\rm M}_{T,\de,\eps,r}\cap(\Ai_m\times\Ai_m)\neq\emptyset\big].
\ec
Inserting this into the left-hand side of (\ref{Plimcros}) and using (\ref{crotig}) we see that the right-hand side of (\ref{Plimcros}) is zero, completing the proof of~I.

It remains to prove II. We will prove that if (\ref{tight})~(i) holds and (\ref{crotig}) fails, then there exist $\mu$ in the closure of $\Mi$ that are not concentrated on $\Ki_{\rm nc}(\Pi^|)$. Fix $\de_n>0$ such that $\de_n\to 0$. Since (\ref{crotig}) does not hold, there exist $T,\eps,\eta>0$ and $r\in\R$ such that for each $n\geq 1$ we can find a random variable $\Ai_n$ with law $\mu_n\in\Mi$ such that
\be\label{mis}
\P\big[\Ci^{\rm M}_{T,\de_n,\eps,r}\cap(\Ai_n\times\Ai_n)\neq\emptyset\big]\geq\eta.
\ee
Since (\ref{tight})~(i) holds, by Theorem~\ref{T:tight}, by going to a subsequence if necessary, we can assume that the $\mu_n$ converge weakly to some probability law $\mu$ on $\Ki_+(\Pi)$. Let $\Ai$ have law $\mu$. By Skorokhod's representation theorem, we can couple $\Ai_n,\Ai$ in such a way that $\Ai_n\to\Ai$ a.s. Then (\ref{mis}) implies that
\be
\P\big[\forall m\in\N\ \exists n\geq m\mbox{ s.t.\ }\Ci^{\rm M}_{T,\de_n,\eps,r}\cap(\Ai_n\times\Ai_n)\neq\emptyset\big]\geq\eta,
\ee
and hence
\be
\P\big[\exists T,\eps>0,\ r\in\R\mbox{ s.t.\ } \forall\de>0,\ m\in\N\ \exists n\geq m\mbox{ s.t.\ }\Ci^{\rm M}_{T,\de,\eps,r}\cap(\Ai_n\times\Ai_n)\neq\emptyset\big]\geq\eta.
\ee
By condition (\ref{limnoncros}) of Lemma~\ref{L:limcros}, this shows that
\be
\P\big[\Ai\mbox{ is non-crossing}]\leq 1-\eta,
\ee
so $\mu$ is not concentrated on $\Ki_{\rm nc}(\Pi^|)$.
\epro

\bpro[of Lemma~\ref{L:biclos}]
Assume that $\Ai_n\in\Ki_{\rm nc}(\Pi^\updo)$ and that $\Ai_n\to\Ai$ for some $\Ai\in\Ki_+(\Pi)$. Then by Lemma~\ref{L:Hauconv}, for each $\pi\in\Ai$, there exist $\pi_n\in\Ai_n$ such that $\pi_n\to\pi$. It is easy to see that $\Pi^\updo$ is a closed subset of $\Pi$, so $\pi\in\Pi^\updo$ for all $\pi\in\Ai$. It remains to show that $\Ai$ is non-crossing. Let $\pi_1,\pi_2\in\Ai$. By Lemma~\ref{L:Hauconv}, there exist $\pi_1^n,\pi_2^n\in\Ai_n$ such that $\pi_i^n\to\pi_i$ $(i=1,2)$. Since $\Ai_n$ is non-crossing, for each $n$, the path $\pi_1^n$ does not cross $\pi_2^n$. By Proposition~\ref{P:collis}, it follows that $\pi_1$ does not cross $\pi_2$, proving that $\Ai$ is non-crossing.
\epro

\end{document}